\documentclass[12pt]{article}
\usepackage{amssymb,amsmath,latexsym}
\usepackage{url}

\title{Weak Convergence of the Scaled Median\\
  of Independent Brownian Motions
  \thanks{This work was supported in part by the VIGRE grants of
    both University of Washington and University of
    Wisconsin-Madison.}}
\author{Jason Swanson\\
  Mathematics Department\\
  University of Wisconsin-Madison\\
  \url{http://www.math.wisc.edu/~swanson}}
\date{June 25, 2006}

\begin{document}

\newtheorem{thm}{Theorem}[section]
\newtheorem{corollary}[thm]{Corollary}
\newtheorem{prop}[thm]{Proposition}
\newtheorem{lemma}[thm]{Lemma}
\newtheorem{remark}{Remark}[section]
\newtheorem{defn}{Definition}[section]

\numberwithin{equation}{section}

\def\MM{\mathcal{M}}
\def\SS{\mathcal{S}}
\def\II{\mathcal{I}}

\def\RR{\mathbb{R}}
\def\NN{\mathbb{N}}
\def\ZZ{\mathbb{Z}}

\def\al{\alpha}
\def\del{\delta}
\def\Del{\Delta}
\def\ee{\varepsilon}
\def\gam{\gamma}
\def\kap{\kappa}
\def\om{\omega}
\def\ph{\varphi}
\def\sig{\sigma}
\def\th{\theta}

\def\pa{\partial}

\def\pf{\noindent{\bf Proof:} }
\def\qed{{$\Box$ \bigskip}}

\def\pt{{\tilde p}}
\def\eet{{\tilde\ee}}
\def\mut{{\tilde\mu}}
\def\At{{\tilde A}}
\def\Ct{{\tilde C}}
\def\Yt{{\tilde Y}}
\def\St{{\tilde S}}
\def\Rt{{\tilde R}}

\providecommand{\flr}[1]{\lfloor#1\rfloor}

\maketitle

\begin{abstract}
We consider the median of $n$ independent Brownian motions,
denoted by $M_n(t)$, and show that $\sqrt{n}\,M_n$ converges
weakly to a centered Gaussian process. The chief difficulty is
establishing tightness, which is proved through direct estimates
on the increments of the median process. An explicit formula is
given for the covariance function of the limit process. The limit
process is also shown to be H\"older continuous with exponent
$\gam$ for all $\gam<1/4$.



\end{abstract}

\section{Introduction}

Consider a dye diffusing in a homogeneous medium. When we view
this phenomenon from a macroscopic perspective, what we see is a
deterministic evolution of the density of the dye, governed by a
partial differential equation. It is well understood that the
solution of this equation can be represented probabilistically in
terms of Brownian motion. The reason, of course, that Brownian
motion enters into this situation is that, heuristically, we can
imagine that each dye particle is performing such a random motion.
In reality, however, a more accurate description of the particles
is that they are following piece-wise linear trajectories and
interacting through collisions.

In 1968, F. Spitzer \cite{S} provided a rigorous connection
between a certain colliding particle model and the Brownian motion
heuristics. In Spitzer's model, we begin with countably many
particles distributed along the real line according to a Poisson
distribution. At time $t=0$, the particles begin moving with
random velocities. These velocities are i.i.d., integrable, mean
zero random variables. During their motion, the particles interact
through elastic collisions. That is, whenever two particles meet,
they exchange velocities (or, equivalently, they exchange
trajectories). The particle which is closest to the origin at time
$t=0$ is called the ``tagged" particle and we denote its position
at time $t$ by $X(t)$. Spitzer showed that the law on
$C[0,\infty)$ induced by the process $t\mapsto c^{-1/2}X(ct)$
converges weakly as $c\to\infty$ to the law of Brownian motion.

Spitzer's work was preceded by that of T. E. Harris \cite{H} who
showed that if the underlying motion of the particles is Brownian,
instead of linear, then $c^{-1/4}X(ct)$ converges to fractional
Brownian motion with Hurst parameter $H=1/4$. These results were
further generalized by D\"urr, Goldstein, and Lebowitz \cite{D} in
1985. They showed, among other things, that if the individual
particles perform fractional Brownian motion with Hurst parameter
$H$, then $c^{-H/2}X(ct)$ converges to fractional Brownian motion
with Hurst parameter $H/2$.

One thing to note in these more general models is the definition
of an ``elastic collision." When the particles perform Brownian
motion, for example, the collisions are not isolated and it is not
entirely clear how to exchange their trajectories at each point of
intersection. In these situations, we generate the collision
process by simply relabelling the particles at each time $t$ in a
way that preserves their initial ordering. For instance, if there
are only finitely particles, as there will be in our model, the
location of the tagged particle is simply an order statistic of
the locations of all of the particles. (In our model, it will be
the median.)

In the work of Spitzer, Harris, and D\"urr et al, the chief
difficulty in proving convergence is establishing tightness. And
in each of these models, the Poisson distribution of the initial
particle configuration provides for tractable computations and is
a central feature of the proofs. In this article, we will consider
a model similar to Harris's, but without the initial Poisson
distribution. Namely, we consider a sequence $\{B_j\}$ of
independent Brownian motions starting at the origin. We let $M_n$
denote the median of the first $n$ of these, and study the scaled
process $X_n=\sqrt{n}\,M_n$. As with the other models, our chief
difficulty will be to prove tightness. We will prove this,
however, by making direct estimates on the path of the ``tagged"
particle, without relying on any special features of the initial
particle distribution. In the end, we will discover a limit
process which behaves locally like fractional Brownian motion with
Hurst parameter $H=1/4$. This fact, formally stated in Theorem
\ref{T:main}, lends support to the evident notion that Harris's
initial Poisson distribution is, to a certain degree, just a
technical convenience, and does not play a significant role in
determining the local behavior of the limit.

Before proceeding with the formal analysis of our model, let us
preview some of the techniques in the proof. The first key
ingredient in the proof will be given by Theorem \ref{T:cond},
which establishes a formula for the conditional law of the median
in terms of probabilities associated with a certain random walk.
The second ingredient will be Lemma \ref{L:chebyplus}, which gives
estimates for this random walk in terms of its parameters. And the
third ingredient will be Lemma \ref{L:taylor} (and its
modification, Lemma \ref{L:taylormod}) which estimates those
parameters in terms of the motion of the individual particles.

Since it would be natural to conjecture that the results of
Spitzer and D\"urr et al would also hold in more general models,
it is important to try to understand how these techniques might
apply in a broader context. For example, we could try to
generalize the results of D\"urr et al by replacing the Brownian
motions in our model with fractional Brownian motions. Or we could
replace them with reflected processes if we wanted to consider
particles in a ``box," reflecting off the walls of the box as well
as each other. Such a model (in which the particles' paths are
piece-wise linear) was studied by P. F. Tupper in \cite{T},
although in that paper, a seemingly ad-hoc condition is imposed in
order to prove tightness. (See the discussion after Theorem 2.3 in
\cite{T}.) Other ways to generalize the model include giving our
particles some nontrivial initial distribution, instead of
starting them at the origin, or possibly considering a quantile
(or even a family of quantiles) other than the median.

In any of these generalized models, the first and second
ingredients outlined above would likely carry over with at most
minor modifications. It is the third ingredient that would not
transfer so easily. The estimates in Lemma \ref{L:taylor} rely
heavily on the fact that the individual particles are performing
Brownian motion. Conceivably, analogous estimates could be worked
out on a case-by-case basis for each different model under
consideration. But the work of Harris, Spitzer, and D\"urr et al
suggests a deeper connection between the motion of the individual
particles and the limit process. It is my belief that this
connection would make itself known through these estimates.

But whether estimates can found in some general form or must be
developed for each model individually, it is my hope that the
techniques developed here can be used to extend the current family
of results to a much broader range of colliding particle models.

\section{The Model and Main Result}

In our model, we will consider a sequence of independent,
standard, one-dimensional Brownian motions,
$\{B_j(t)\}_{j=1}^\infty$. Let $M_n(t)$ denote the median of the
first $n$ Brownian motions. To be precise, define the median
function $\MM_n:\RR^n\to\RR$ as follows: if
$(x_1,\ldots,x_n)\in\RR^n$ and $\tau$ is a permutation of
$\{1,\ldots,n\}$ such that $x_{\tau(1)}\le x_{\tau(2)}\le
\cdots\le x_{\tau(n)}$, then $\MM_n(x_1,\ldots,x_n)=x_{\tau(k)}$,
where $k=\flr{(n+1)/2}$ and $\flr{x}$ denotes the greatest integer
less than or equal to $x$. We then define the (continuous) median
process $M_n(t)=\MM_n(B_1(t),\ldots,B_n(t))$.

In terms of colliding particles, what we have here is a sequence
of particle systems. In the $n$-th system there are $n$ particles
performing Brownian motion. If these particles interact through
elastic collisions, then their trajectories are given by the order
statistics of $B_1(t),\ldots,B_n(t)$. We will investigate the
behavior of the center particle's trajectory, $M_n(t)$.

In order to get a non-degenerate limit, we must consider the
scaled median process $X_n(t)=\sqrt{n}\,M_n(t)$. The random
variables $X_n=\{X_n(t):0\le t<\infty\}$ take values in the space
$C[0,\infty)$, which we endow with the topology of uniform
convergence on compact sets. It will be shown that these processes
converge weakly, by which we mean that they converge in law as
$C[0,\infty)$-valued random variables.

\begin{thm}\label{T:main}
There exists a continuous process $X=\{X(t):0\le t<\infty\}$ such
that $X_{2n+1}$ converges weakly to $X$ as $n\to\infty$. Moreover,
$X$ is a centered Gaussian process, which is locally H\"older
continuous with exponent $\gam$ for every $\gam\in(0,1/4)$, and
has covariance function
  \begin{equation}
  E[X(s)X(t)]=\sqrt{st}\,\sin^{-1}%
    \left({\frac{s\wedge t}{\sqrt{st}}}\right),\label{8a}
  \end{equation}
where $\sin^{-1}(\cdot)$ takes values in $[-\pi/2,\pi/2]$.
\end{thm}

It can be shown by \eqref{8a} that, for $t-s$ small,
$E|X(t)-X(s)|^2\approx\sqrt{t-s}$. In other words, the limit
process has the same local fluctuations as fractional Brownian
motion with Hurst parameter $H=1/4$.

The chief difficulty in proving Theorem \ref{T:main} will be to
establish the tightness of the processes $X_{2n+1}$. Before
dealing with this issue, let us first establish the convergence of
the finite-dimensional distributions and the existence of the
limit process. To begin, we will need the following result, which
is a special case of Theorems 7.1.1 and 7.1.2 in \cite{R}.

\begin{thm}\label{T:MCLT}
Let $\{\xi^{(n)}\}_{n=1}^\infty$ be an i.i.d. sequence of random
vectors in $\RR^d$ and define the component-wise median of
$\xi^{(1)},\ldots,\xi^{(n)}$ to be the vector $M^{(n)}$ with
components
$M^{(n)}_j=\MM_n(\xi^{(1)}_j,\xi^{(2)}_j,\ldots,\xi^{(n)}_j)$. Let
$F_j(x)=P(\xi^{(1)}_j\le x)$, %
$G_{ij}(x,y)=P(\xi^{(1)}_i\le x,\xi^{(1)}_j\le y)$, and
$\rho_{ij}=G_{ij}(0,0)-1/4$. If
  \begin{description}
  \item{(i)} $F_j(0)=1/2$ and $F_j'(0)>0$ for all $j$, and%
  \item{(ii)} $G_{ij}$ is continuous at $(0,0)$ for all $i$ and $j$,
  \end{description}
then $\sqrt{n}\,M^{(n)}$ converges in law to a jointly Gaussian
random vector $N$ satisfying
  \[
  EN_iN_j=\frac{\rho_{ij}}{F_i'(0)F_j'(0)}
  \]
and $EN_i=0$.
\end{thm}

For our purposes, we will need the following.

\begin{corollary}\label{C:MCLT}
If $\{\xi^{(n)}\}_{n=1}^\infty$ is an i.i.d. sequence of jointly
Gaussian random vectors in $\RR^d$ with mean zero and covariance
matrix $\sig$, then $\sqrt{n}\,M^{(n)}$ converges in law to a
jointly Gaussian random vector $Z$ with mean zero and covariance
matrix $\tau$, where
  \[
  \tau_{ij}=EZ_iZ_j=\sqrt{\sig_{ii}\sig_{jj}}\,\sin^{-1}
    \left({\frac{\sig_{ij}}{\sqrt{\sig_{ii}\sig_{jj}}}}\right)
  \]
and $\sin^{-1}(\cdot)$ takes values in $[-\pi/2,\pi/2]$.
\end{corollary}

\pf This follows easily from Theorem \ref{T:MCLT} and the
well-known fact that if $X$ and $Y$ are jointly Gaussian with mean
zero, then
  \[
  P(X\le0,Y\le0)=\frac14+\frac1{2\pi}\sin^{-1}
    \left({\frac{EXY}{\sqrt{EX^2\cdot EY^2}}}\right),
  \]
where $\sin^{-1}(\cdot)$ takes values in $[-\pi/2,\pi/2]$. \qed

\begin{thm}\label{T:limit}
There exists a centered Gaussian process %
$X=\{X(t):0\le t<\infty\}$ with covariance function \eqref{8a} and
which is locally H\"older continuous with exponent $\gam$ for
every $\gam\in(0,1/4)$.
\end{thm}

\pf Let $T$ be the set of finite sequences %
${\bf t}=(t_1,\ldots,t_n)$ of distinct, nonnegative numbers, where
the length $n$ of these sequences ranges over the set of positive
integers. For each ${\bf t}$ of length $n$, let %
${\bf Z}_{\bf t}=(Z_1,\ldots,Z_n)$ be a jointly Gaussian random
vector with mean zero and covariance
  \[
  EZ_iZ_j=\sqrt{t_it_j}\,\sin^{-1}%
    \left({\frac{t_i\wedge t_j}{\sqrt{t_it_j}}}\right).
  \]
(By Corollary \ref{C:MCLT}, with
$\xi^{(j)}=(B_j(t_1),\ldots,B_j(t_n))$, such a ${\bf Z}_{\bf t}$
exists.) Define the measure $Q_{\bf t}$ on $\RR^n$ by %
$Q_{\bf t}(A)=P({\bf Z}_{\bf t}\in A)$. The family of
finite-dimensional distributions, $\{Q_{\bf t}\}_{{\bf t}\in T}$,
is clearly consistent, so there exists a real-valued process
$X=\{X(t):0\le t<\infty\}$ that has the desired finite-dimensional
distributions. It remains only to show that this process has a
continuous modification, which is locally H\"older-continuous with
exponent $\gam$ for every $\gam\in(0,1/4)$.

By the Kolmogorov-\v{C}entsov Theorem (Theorem 2.2.8 in \cite{K}),
if, for each $T>0$,
  \[
  E|X(t)-X(s)|^\al\le C_T|t-s|^{1+\beta}
  \]
for some positive constants $\al$, $\beta$, and $C_T$ (depending
on $T$) and all $0\le s<t\le T$, then $X$ has a continuous
modification which is locally H\"older-continuous with exponent
$\gam$ for every $\gam\in(0,\beta/\al)$. Hence, it will suffice
for us to show that for every $\al>4$ and every $T>0$,
  \[
  E|X(t)-X(s)|^\al\le C|t-s|^{\al/4}
  \]
for some $C>0$ (depending only on $T$ and $\al$) and all %
$0\le s<t\le T$.

First, observe that $X(t)-X(s)$ is normal with mean zero and
variance
  \begin{align*}
  E|X(t)-X(s)|^2 &=EX(t)^2+EX(s)^2-2EX(t)X(s)\\
  &=\frac{\pi}2t+\frac{\pi}2s
    -2\sqrt{st}\,\sin^{-1}\left({\sqrt{\frac st}}\right).
  \end{align*}
An application of L'H\^{o}pital's Rule shows that
  \[
  \frac{\pi/2-\sin^{-1}x}{\sqrt{1-x^2}}\to1
  \]
as $x\to 1$. Hence, for some positive constant $C'$, we have
$-\sin^{-1}x\le C'\sqrt{1-x^2}-\pi/2$ for all $0\le x\le 1$. Now
let $x=s/t$. Then
  \begin{align*}
  E|X(t)-X(s)|^2 &=t\left[{
    \frac{\pi}2+\frac{\pi}2x-2\sqrt{x}\,\sin^{-1}(\sqrt{x})%
    }\right]\\
  &\le t\left[{\frac{\pi}2+\frac{\pi}2x+2\sqrt{x}\left({
    C'\sqrt{1-x}-\frac{\pi}2}\right)}\right]\\
  &=t\left[{\frac{\pi}2(1-\sqrt{x})^2+2C'\sqrt{x}\sqrt{1-x}
    }\right].
  \end{align*}
Since $1-\sqrt{x}\le\sqrt{1-x}$ for $0\le{x}\le1$,
  \begin{align*}
  E|X(t)-X(s)|^2 &\le t\left({
    \frac{\pi}2(1-x)+2C'\sqrt{x}\sqrt{1-x}}\right)\\
  &\le t\left({\frac{\pi}2\sqrt{1-x}+2C'\sqrt{1-x}}\right)\\
  &=\sqrt{t}\left({\frac{\pi}2+2C'}\right)\sqrt{t-s}\\
  &\le C''|t-s|^{1/2},
  \end{align*}
where $C''=\sqrt{T}(\pi/2+2C')$.

Now, for every $\al>0$, there is a constant $K_\al$ such that if
$N$ is normal with $EN=0$, then $E|N|^\al=K_\al(EN^2)^{\al/2}$.
Thus, for any $\al>4$, $E|X(t)-X(s)|^\al\le C|t-s|^{\al/4}$, where
$C=K_\al(C'')^{\al/2}$. \qed

\begin{thm}\label{T:fdd}
Let $X(t)$ be as in Theorem \ref{T:limit} and let %
$0\le t_1<\cdots<t_d$, $d\ge1$, be arbitrary. Then
$(X_n(t_1),\ldots,X_n(t_d))$ converges in law to
$(X(t_1),\ldots,X(t_d))$ as $n\to\infty$.
\end{thm}

\pf This is an immediate consequence of Corollary \ref{C:MCLT}.
\qed

It now follows (see, for example, Theorem 2.4.15 in \cite{K}) that
Theorem \ref{T:main} will be proved once we establish the
following result.

\begin{thm}\label{T:tight}
The sequence of processes $\{X_{2n+1}\}_{n=1}^\infty$ is tight.
\end{thm}

\section{Conditions for Tightness}

A sufficient condition for tightness which will serve as the
starting point for our analysis is the following.

\begin{thm}\label{T:moment}
If $\{Z_n\}$ is a sequence of continuous stochastic processes such
that
  \begin{description}
  \item{(i)} $\sup_nP(|Z_n(t)-Z_n(s)|\ge\ee)%
    \le C_T\ee^{-\al}|t-s|^{1+\beta}$ %
    whenever $0<\ee<1$, $T>0$, and $0\le s,t\le T$, and%
  \item{(ii)} $\sup_nE|Z_n(0)|^\nu<\infty$
  \end{description}
for some positive constants $\al$, $\beta$, $\nu$, and $C_T$
(depending on $T$), then $\{Z_n\}$ is tight.
\end{thm}

An alternative formulation of this theorem is one in which
condition (i) is replaced by
  \begin{equation}
  \sup_{n\ge1}E|Z_n(t)-Z_n(s)|^\al%
    \le C_T|t-s|^{1+\beta}.\label{16a}
  \end{equation}
For a proof of this alternative version, the reader is referred to
Problem 2.4.11 in \cite{K}, which has a worked solution. An
inspection of the proof shows that \eqref{16a} is needed only to
establish (via Chebyshev's inequality) condition (i).

Since the median process inherits the scaling property of Brownian
motion, we will find it convenient to reformulate Theorem
\ref{T:moment}. Specifically, for any real number $c\ge0$ and any
$x\in\RR^d$, we have $\MM_n(cx)=c\MM_n(x)$. Hence, the processes
$X_n(c\,\cdot)$ and $\sqrt{c}\,X_n(\cdot)$ have the same law. For
processes with this scaling property, we can modify Theorem
\ref{T:moment} in the following way.

\begin{thm}\label{T:moment2}
Let $\{Z_n\}$ be a sequence of continuous stochastic processes.
Suppose there exists $r>0$ such that for every $c\ge0$ and every
$n$, the processes $Z_n(c\,\cdot)$ and $c^rZ_n(\cdot)$ have the
same law. Suppose further that
  \begin{description}
  \item{(i)} $\sup_nP(|Z_n(1+\del)-Z_n(1)|>\ee)%
    \le C\ee^{-\al}\del^{1+\beta}$ %
    whenever $0<\ee<1$ and $0<\del<\del_0$
  \end{description}
for some positive constants $\del_0$, $C$, $\al$, and $\beta$.
Define $\gam=\min(\al r,\beta r,1+\beta)$. If $\gam>1$ and
  \begin{description}
  \item{(ii)} $\sup_nE|Z_n(1)|^{\gam/r}<\infty$,
  \end{description}
then the sequence $\{Z_n\}$ is tight.
\end{thm}

Theorem \ref{T:moment2} follows directly from Theorem
\ref{T:moment} (a complete proof can be found starting on p.36 of
\cite{Sw}). We will be applying it to the sequence $Z_n=X_{2n+1}$,
in which case we have $r=1/2$. We will find it quite
straightforward to verify condition (ii). To verify condition (i),
we will utilize the following lemma, which will be the central
focus of the remainder of our analysis.

\begin{lemma}\label{L:key}
There exists a constant $\del_0>0$ and a family of constants
$\{C_p\}_{p>2}$ such that for each $p>2$,
  \begin{equation}
  \sup_{n\ge3}P\left({M_n(1+\del)-M_n(1)%
    >\frac\ee{\sqrt{n}}}\right)%
    \le C_p(\ee^{-1}\del^{1/6})^p\label{21a}
  \end{equation}
whenever $0<\ee<1$ and $0<\del\le\del_0$.
\end{lemma}

It has already been remarked that the limit process $X$ behaves
locally like a fractional Brownian motion with Hurst parameter
$H=1/4$. It seems reasonable, then, to conjecture that the
right-hand side of \eqref{21a} could be replaced by
$C_p(\ee^{-1}\del^{1/4})^p$. Although this sharper bound was not
obtained, the choice of $1/6$ as the exponent in \eqref{21a}
appears to be arbitrary. Presumably, with minor modifications to
the proofs presented here, the right-hand side of \eqref{21a}
could be replaced by $C_p(\ee^{-1}\del^\nu)^p$ for any fixed
$\nu<1/4$.

\bigskip\noindent{\bf Proof of Theorem \ref{T:tight}, given Lemma
\ref{L:key}. } We apply Theorem \ref{T:moment2} to $Z_n=X_{2n+1}$
with $r=1/2$. Choose any $p>18$, let $\al=p$, and let
$\beta=(p-6)/6$. Note that, in this case, $\gam=\beta/2>1$.

To verify condition (i), let $\del_0$ be as in Lemma \ref{L:key}.
Since $X_{2n+1}(\cdot)$ and $-X_{2n+1}(\cdot)$ have the same law,
  \begin{align*}
  \sup_{n\ge1}P(|X_{2n+1}(1+\del)-X_{2n+1}(1)|>\ee)
    &=2\sup_{n\ge1}P(X_{2n+1}(1+\del)-X_{2n+1}(1)>\ee)\\
  &\le2C_p(\ee^{-1}\del^{1/6})^p\\
  &=2C_p\ee^{-\al}\del^{1+\beta}
  \end{align*}
whenever $0<\ee<1$ and $0<\del<\del_0$.

To verify condition (ii), we will show that for any $q>0$,
  \[
  \sup_{n\ge1}E|X_{2n+1}(1)|^q<\infty.
  \]
To see this, observe that for $n$ odd,
  \begin{align*}
  E|X_n(1)|^q &=\int_0^\infty{qy^{q-1}P(|X_n(1)|>y)\,dy}\\
  &=2\int_0^\infty{qy^{q-1}P(X_n(1)<-y)\,dy}.
  \end{align*}
It will therefore suffice to show that for any $\kap>2$, there
exists a finite constant $K$ such that
  \begin{equation}
  P(X_n(1)<-y)\le Ky^{-\kap}\label{24a}
  \end{equation}
for all $y>0$ and all $n$.

To prove \eqref{24a}, we will consider two cases. First, assume
$y\ge2\sqrt{n}$. Note that by Theorem 1.3.2 in \cite{R}, $M_n(1)$
has density
  \begin{equation}
  f_n(x)=k\binom{n}{k}\frac1{2\pi}\,
    \Phi(x)^{k-1}\Phi(-x)^{n-k}e^{-x^2/2}\label{25a}
  \end{equation}
where $k=\flr{(n+1)/2}$ and
$\Phi(x)=\frac1{2\pi}\int_{-\infty}^x{e^{-u^2/2}\,du}$. Hence,
  \begin{align*}
  P(X_n(1)<-y) &=P(M_n(1)<-y/\sqrt{n})\\
  &=\frac{n!}{(n-k)!(k-1)!}\int_{-\infty}^{-y/\sqrt{n}}
    {\Phi(x)^{k-1}\Phi(-x)^{n-k}\Phi'(x)\,dx}\\
  &\le\frac{n^k}{(k-1)!}\int_{-\infty}^{-y/\sqrt{n}}
    {\Phi(x)^{k-1}\Phi'(x)\,dx}\\
  &=\frac{n^k}{k!}\Phi(-y/\sqrt{n})^k.
  \end{align*}
By Stirling's formula, there exists a universal positive constant
$C$ such that $k!\ge C^{-1}k^ke^{-k}$. Also, writing
$\int_x^\infty{e^{-u^2/2}\,du}%
=\int_x^\infty{u^{-1}\cdot ue^{-u^2/2}\,du}$ and integrating by
parts, it follows that
  \begin{equation}
  \sqrt{2\pi}\,\Phi(-x)\le x^{-1}e^{-x^2/2}\label{25b}
  \end{equation}
for all $x>0$. Thus,
  \[
  P(X_n(1)<-y)\le C\frac{n^k}{k^ke^{-k}}
    \left({\frac{\sqrt{n}}{y}\,e^{-y^2/2n}}\right)^k.
  \]
Since $y\ge2\sqrt{n}$ and $n/k\le2$, we have
  \[
  P(X_n(1)<-y)\le Ce^{k(1-y^2/2n)}.
  \]
Since $1-y^2/2n<0$ and $k\ge n/2$, we have
  \[
  P(X_n(1)<-y)\le Ce^{n/2-y^2/4}\le Ce^{-y^2/8}.
  \]
Finally, given $\kap>2$, there exists $K$ such that
$Ce^{-y^2/8}\le Ky^{-\kap}$ for all $y>0$, which verifies
\eqref{24a} in the case $y\ge2\sqrt{n}$.

Now assume $y<2\sqrt{n}$. In this case,
  \begin{align*}
  P(X_n(1)<-y) &=P(M_n(1)<-y/\sqrt{n})\\
  &=P\biggl({\sum_{j=1}^n
    1_{\{B_j(1)<-y/\sqrt{n}\}}\ge\tfrac n2}\biggr)\\
  &=P\biggl({\sum_{j=1}^n%
    \xi_j\ge n\left({\tfrac12-\mu}\right)}\biggr),
  \end{align*}
where $\mu=\Phi(-y/\sqrt{n})$ and
$\xi_j=1_{\{B_j(1)<-y/\sqrt{n}\}}-\mu$. By Burkholder's inequality
(see, for example, Theorem 6.3.10 in \cite{St}), there exists a
constant $K'$, depending only on $\kap$, such that
  \[
  E\Bigl|{\sum_{j=1}^n\xi_j}\Bigr|^\kap \le
    K'E\Bigl|{\sum_{j=1}^n|\xi_j|^2}\Bigr|^{\kap/2}.
  \]
Hence, since $\kap>2$, Jensen's inequality and the fact that
$|\xi_j|\le1$ a.s. imply
  \[
  E\Bigl|{\sum_{j=1}^n\xi_j}\Bigr|^\kap
    \le K'n^{\kap/2}E\sum_{j=1}^n\frac1n|\xi_j|^k
    \le K'n^{\kap/2}.
  \]
Chebyshev's inequality now gives
  \[
  P(X_n(1)<-y)\le\frac{K'n^{\kap/2}}
    {\left|{n\left({\frac12-\mu}\right)}\right|^\kap}
    =K'\left|{\sqrt{n}\left({\tfrac12-\mu}\right)}\right|^{-\kap}.
  \]
Since
  \[
  \sqrt{n}\left({\tfrac12-\mu}\right)%
    =\frac{\sqrt{n}}{\sqrt{2\pi}}%
    \int_0^{y/\sqrt{n}}{e^{-u^2/2}\,du}%
    \ge\frac y{\sqrt{2\pi}}\,e^{-y^2/2n}%
    \ge\frac y{\sqrt{2\pi}}\,e^{-2},
  \]
we have that $P(X_n(1)<-y)\le Ky^{-\kap}$, where
$K=K'(e^{-2}/\sqrt{2\pi})^{-\kap}$. This verifies \eqref{24a} when
$y<2\sqrt{n}$ and completes the proof. \qed

Our goal for the remainder of this article is to establish
\eqref{21a}. Since each individual Brownian particle can be
expected to move a distance of $\sqrt{\del}$ between time $t=1$
and $t=1+\del$, we will accomplish our goal by considering three
different ``jump regimes." They are: the large jump regime in
which $\ee/\sqrt{n}$ is much larger than $\sqrt{\del}$, the small
jump regime in which $\ee/\sqrt{n}$ is much smaller than
$\sqrt{\del}$, and the medium jump regime in which these two
quantities are comparable. In the first two regimes, we will
establish the sharp bound mentioned in the remark following Lemma
\ref{L:key}. The bound in the medium jump regime will be
established by modifying the techniques used in the small jump
regime. This modification will result in the weaker bound given in
\eqref{21a}.

\section{The Large Jump Regime}

The large jump regime is the easiest of the three to deal with.
The probability that the median makes a large jump can be bounded
above by the probability that at least one Brownian particle makes
a large jump. Since the latter probability is exponentially small,
the derivation of \eqref{21a} is immediate.

\begin{lemma}\label{L:large}
Fix $p>0$ and $0<\Del<1/2$. Suppose that $\ee,\del\in(0,1)$ and
$n\in\NN$ satisfy $\ee/\sqrt{n}\ge\del^{1/2-\Del}$. Then
  \[
  P\left({M_n(1+\del)-M_n(1)>\frac\ee{\sqrt{n}}}\right) \le
    C(\ee^{-1}\del^{1/4})^p,
  \]
where $C$ depends only on $p$ and $\Del$.
\end{lemma}

\pf Suppose that $B_j(1+\del,\om)-B_j(1,\om)\le\ee/\sqrt{n}$ for
all $j$. Then, for each $j$ such that $B_j(1,\om)\le M_n(1,\om)$,
we have $B_j(1+\del,\om)\le M_n(1,\om)+\ee/\sqrt{n}$. Note that
there are at least $k=\flr{(n+1)/2}$ such values of $j$. It
follows that $M_n(1+\del,\om)\le M_n(1,\om)+\ee/\sqrt{n}$.
Therefore,
  \[
  \bigcap_{j=1}^n\{B_j(1+\del)-B_j(1)\le\ee/\sqrt{n}\}
  \subset\{M_n(1+\del)-M_n(1)\le \ee/\sqrt{n}\},
  \]
which gives
  \begin{align*}
  P\left({M_n(1+\del)-M_n(1)>\frac\ee{\sqrt{n}}}\right)%
    &\le P\biggl({\bigcup_{j=1}^n
    \{B_j(1+\del)-B_j(1)>\ee/\sqrt{n}\}}\biggr)\\
  &\le n\Phi(-\ee/\sqrt{n\del}).
  \end{align*}
For each $r>0$, there exists $C_r$ such that%
$\Phi(-x)\le C_rx^{-r}$ for all $x>0$. Taking $r=(p/4+1)/\Del$
gives
  \[
  P\left({M_n(1+\del)-M_n(1)>\frac\ee{\sqrt{n}}}\right)%
    \le nC_r\left({\frac\ee{\sqrt{n\del}}}\right)^{-r}%
    \le nC_r(\del^{-\Del})^{-r}%
    =C_rn\del^{p/4+1}.
  \]
The proof is completed by observing that
$n\le\ee^2\del^{-1}\le\ee^{-p}\del^{-1}$. \qed

This establishes the necessary bound for the large jump regime.
The other regimes, as we will see, are considerably more difficult
to deal with.

\section{Conditioning the Median}\label{cond}

To establish \eqref{21a} for the small and medium jump regimes, we
will use conditioning. It may seem natural, at first, to condition
on the locations of all the Brownian particles at time $t=1$. It
turns out, however, that this is, in some sense, too much
information. Rather, we shall condition only on the location of
the median particle at time $t=1$.

Let us first give a heuristic description of this conditioning.
Suppose that $M_n(1)=x$. This tells us three things. First, we
have a single Brownian particle whose location is $x$. Second, we
have roughly $n/2$ Brownian particles whose locations are less
than $x$. Other than this condition on their locations, these
particles are independent and identically distributed. We will
refer to these particles as the ``lower" particles. Third, we have
roughly $n/2$ i.i.d. Brownian particles whose locations are
greater than $x$. These will naturally be referred to as the
``upper" particles.

Let us now fix $y>0$ and consider the event
$D=\{M_n(1+\del)-M_n(1)>y\}$. This event will occur if and only if
there are at least $n/2$ particles whose location at time
$t=1+\del$ is greater than $x+y$. Particles that satisfy this
condition will be said to have ``jumped." Let $U(j)$ be the event
that the $j$-th upper particle jumps, and let $L(j)$ be the event
that the $j$-th lower particle does {\it not} jump. Then the total
number of particles that jump is
  \[
  \sum1_{U(j)}+\left({\frac n2-\sum1_{L(j)}}\right).
  \]
The event $D$ will occur if and only if this sum is at least
$n/2$, which occurs if and only if $\sum Y_j\ge0$, where
$Y_j=1_{U(j)}-1_{L(j)}$ are i.i.d. $\{-1,0,1\}$-valued random
variables. Through conditioning, then, we are able to transform
the event of interest into one involving an i.i.d. sum.

With these heuristics in place, let us establish the rigorous
result. Define
  \begin{align}
  p_1=p_1(x,y,\del)&=P(B(1+\del)<x+y|B(1)<x)\label{37a}\\
  p_2=p_2(x,y,\del)&=P(B(1+\del)>x+y|B(1)>x)\label{37b}\\
  &=p_1(-x,-y,\del)\notag
  \end{align}
and
  \begin{equation}
  q_j=1-p_j.\label{37c}
  \end{equation}
In the language of our heuristics, $p_1$ is the probability that a
lower particle does not jump and $p_2$ is the probability that an
upper particle does jump.

Now, for each fixed triple $(x,y,\del)$, let
$\{\xi^L_j\}_{j=1}^\infty$ and $\{\xi^U_j\}_{j=1}^\infty$ be
sequences of i.i.d $\{0,1\}$-valued random variables with
$P(\xi^L_j=1)=p_1$ and $P(\xi^U_j=1)=p_2$. Define
$Y_j=\xi^U_j-\xi^L_j$. Observe that $\{Y_j\}_{j=1}^\infty$ is an
i.i.d. sequence of $\{-1,0,1\}$-valued random variables and, for
future reference, define
  \begin{align}
  \pt_1&=P(Y_j=-1)=p_1q_2\label{38a}\\
  \pt_2&=P(Y_j=1)=p_2q_1\label{38b}\\
  \eet&=P(Y_j\ne0)=\pt_1+\pt_2\label{38c}\\
  \mut&=-EY_j=\pt_1-\pt_2.\label{38d}
  \end{align}
Finally, let $S_k=\sum_{j=1}^kY_j$ and
$\ph_k(x,y,\del)=P(S_k\ge0)$.

Our heuristics suggest that
  \[
  P(M_n(1+\del)-M_n(1)>y|M_n(1)=x)\approx\ph_{n/2}(x,y,\del).
  \]
For a rigorous statement, the following inequality will serve our
purposes.

\begin{thm}\label{T:cond}
Let $n\ge3$ and $k=\flr{(n+1)/2}$. Then for all $y>0$ and all
$\del>0$,
  \[
  P(M_n(1+\del)-M_n(1)>y)
    \le\int_{-\infty}^\infty{\ph_{k-1}(x,y,\del)f_n(x)\,dx},
  \]
where $f_n(x)$ is the density of $M_n(1)$, given by \eqref{25a}.
\end{thm}

\pf First, let us observe that
  \begin{align*}
  \ph_k(x,y,\del) =P(S_k\ge0)
    &=P\left({\sum_{j=1}^k\xi^U_j\ge\sum_{j=1}^k\xi^L_j}\right)\\
  &=\sum_{\ell=0}^k\sum_{m=\ell}^k
    P\left({\sum_{j=1}^k\xi^L_j=\ell,\,
    \sum_{j=1}^k\xi^U_j=m}\right)\\
  &=\sum_{\ell=0}^k\sum_{m=\ell}^k\binom{k}{\ell}\binom{k}{m}
    p_1^\ell q_1^{k-\ell}p_2^mq_2^{k-m}.
  \end{align*}
Let us also adopt the following notation: for $h>0$, let
$p_{1,h}=p_1(x+h,y-h,\del)$ and
  \begin{equation*}
  \ph_k^h(x,y,\del)
    =\sum_{\ell=0}^k\sum_{m=\ell}^k\binom{k}{\ell}\binom{k}{m}
    p_{1,h}^\ell q_{1,h}^{k-\ell}p_2^mq_2^{k-m},
  \end{equation*}
where $q_{1,h}=1-p_{1,h}$. Finally, let %
$\Del M_n=M_n(1+\del)-M_n(1)$.

Now, fix $\del>0$ and $y>0$. Let $K\in\NN$ and let $h>0$ with
$K/h\in\NN$. Then
  \[
  P(\Del M_n>y,\,|M_n(1)|\le K)%
    \le\sum_{\substack{x\in h\ZZ\\
    |x|\le K}}%
    P\left({M_n(1+\del)>x+y,\,M_n(1)\in[x,x+h)}\right).
  \]
Let $\SS_n=\{1,\ldots,n\}$ and let $S=S_n$ denote the collection
of all ordered pairs $(I,j)$ where $I\subset\SS_n$ and $j\in\SS_n$
satisfy $|I|=k-1$ and $j\notin I$. For $(I,j)\in S$, $x\in\RR$,
and $h>0$, define $I(j)^c=\SS_n\setminus(I\cup\{j\})$ and
  \begin{equation*}\begin{split}
  A(I,j,x,h)&=\{B_j(1)\in[x,x+h)\}\\
  &\quad\cap\{B_i(1)<B_j(1),\,\forall i\in I\}
    \cap\{B_i(1)>B_j(1),\,\forall i\in I(j)^c\},\\
  \At(I,j,x,h)&=\{B_j(1)\in[x,x+h)\}\\
  &\quad\cap\{B_i(1)<x+h,\,\forall i\in I\}
    \cap\{B_i(1)>x,\,\forall i\in I(j)^c\}.
  \end{split}\end{equation*}
Note that $\{M_n(1)\in[x,x+h)\}=\bigcup\{A(I,j,x,h):(I,j)\in S\}$
up to a set of measure zero, and that this is a disjoint union.
Therefore,
  \begin{align*}
  P(M_n(1+\del)>x+y,\,M_n(1)\in[x,x+h))%
    &=\sum_{(I,j)\in S}%
    P(M_n(1+\del)>x+y,\,A(I,j,x,h))\\
  &\le\sum_{(I,j)\in S}%
    P(M_n(1+\del)>x+y,\,\At(I,j,x,h)),
  \end{align*}
since $A(I,j,x,h)\subset\At(I,j,x,h)$.

Now fix $(I,j)\in S$ and $x\in\RR$. Define
  \begin{align*}
  N_1&=\sum_{i\in I}\,\,1_{\{B_i(1+\del)<x+y\}}\\
  N_2&=\sum_{i\in I(j)^c}1_{\{B_i(1+\del)>x+y\}}\\
  N&=\sum_{i=1}^n\,\,1_{\{B_i(1+\del)>x+y\}}
  \end{align*}
and note that $\{M_n(1+\del)>x+y\}=\{N\ge n-k+1\}$. Also note
that, up to a set of measure zero,
  \begin{align*}
  N&=N_2+(k-1)-N_1+1_{\{B_j(1+\del)>x+y\}}\\
  &\le N_2-N_1+k.
  \end{align*}
Thus, if $d(n)=n-2k+1$, then
$\{M_n(1+\del)>x+y\}\subset\{N_2-N_1\ge d(n)\}$. This gives
  \begin{align*}
  P(M_n(1+\del)>x+y,\,\At(I,j,x,h))%
    &\le P(N_2-N_1\ge d(n),\,\At(I,j,x,h))\\
  &=\sum_{\ell=0}^{k-1}\sum_{m=d(n)+\ell}^{n-k}%
    P(N_1=\ell,\,N_2=m,\,\At(I,j,x,h)).
  \end{align*}
Hence, if we define
  \begin{align*}
  P_1(\ell)&=P(\{N_1=\ell\}\cap\{B_i(1)<x+h,\,\forall i\in I\}),\\
  P_2(m)&=P(\{N_2=m\}\cap\{B_i(1)>x,\,\forall i\in I(j)^c\}),
  \end{align*}
then we can write
  \[
  P(M_n(1+\del)>x+y,\,\At(I,j,x,h))%
    \le\sum_{\ell=0}^{k-1}\sum_{m=d(n)+\ell}^{n-k}%
    P(B_j(1)\in[x,x+h))P_1(\ell)P_2(m).
  \]
Since
  \[
  P(\At(I,j,x,h))%
    =P(B_j(1)\in[x,x+h))\Phi(x+h)^{k-1}\Phi(-x)^{n-k},
  \]
this gives
  \[
  P(M_n(1+\del)>x+y|\At(I,j,x,h))%
    \le\sum_{\ell=0}^{k-1}\sum_{m=d(n)+\ell}^{n-k}%
    \frac{P_1(\ell)}{\Phi(x+h)^{k-1}}\cdot%
    \frac{P_2(m)}{\Phi(-x)^{n-k}}
  \]
for each fixed $I$, $j$, and $x$.

To simplify this double sum, let
  \begin{equation}\label{45a}
  \begin{split}
  \psi(x,y,\del)&=P(B(1+\del)<x+y,\,B(1)<x)\\
  &=\int_{-\infty}^x{\Phi\left({\frac{x+y-t}{\sqrt{\del}}}\right)
    \Phi'(t)\,dt}.
  \end{split}
  \end{equation}
Then by symmetry and independence,
  \begin{align*}
  P_1(\ell)&=\binom{k-1}{\ell}%
    [\psi(x+h,y-h)]^\ell[\Phi(x+h)-\psi(x+h,y-h)]^{k-1-\ell},\\
  P_2(m)&=\binom{n-k}{m}%
    [\psi(-x,-y)]^m[\Phi(-x)-\psi(-x,-y)]^{n-k-m}.
  \end{align*}
Also note that
  \[
  \frac{\psi(x+h,y-h)}{\Phi(x+h)}%
    =P(B(1+\del)<x+y|B(1)<x+h)=p_{1,h}
  \]
and
  \[
  \frac{\psi(-x,-y)}{\Phi(-x)}=P(B(1+\del)>x+y|B(1)>x)=p_2,
  \]
which yields
  \[
  P(M_n(1+\del)>x+y|\At(I,j,x,h))%
    \le\sum_{\ell=0}^{k-1}\sum_{m=d(n)+\ell}^{n-k}%
    \binom{k-1}{\ell}\binom{n-k}{m}%
    p_{1,h}^\ell q_{1,h}^{k-1-\ell}p_2^mq_2^{n-k-m}
  \]
for each fixed $I$, $j$, and $x$.

Now suppose $n$ is odd. In this case, $d(n)=0$ and $n-k=k-1$, so
  \begin{equation}\label{42a}
  P(M_n(1+\del)>x+y|\At(I,j,x,h))\le\ph_{k-1}^h(x,y,\del).
  \end{equation}
On the other hand, if $n$ is even, then $d(n)=1$ and $n-k=k$, so
  \begin{align*}
  P(M_n(1+\del)>x+y|\At(I,j,x,h))%
    &\le\sum_{\ell=0}^{k-1}\sum_{m=\ell+1}^{k} \binom{k-1}{\ell}%
    \binom{k}{m}p_{1,h}^\ell q_{1,h}^{k-1-\ell}p_2^mq_2^{k-m}\\
  &=\sum_{\ell=0}^{k-1}\binom{k-1}{\ell}%
    p_{1,h}^\ell q_{1,h}^{k-1-\ell}%
    \sum_{m=\ell+1}^{k}\binom{k}{m}p_2^mq_2^{k-m}.
  \end{align*}
But
  \[
  \sum_{m=\ell+1}^k\binom{k}{m}p_2^mq_2^{k-m}%
    =P\biggl({\sum_{j=1}^k\xi^U_j>\ell}\biggr)%
    \le P\biggl({\sum_{j=1}^{k-1}\xi^U_j\ge\ell}\biggr)%
    =\sum_{m=\ell}^{k-1}\binom{k-1}{m}p_2^mq_2^{k-1-m},
  \]
so \eqref{42a} holds in this case as well.

Putting it all together, we have
  \begin{align*}
  P(\Del M_n>y,\,|M_n(1)|\le K)%
    &\le\sum_{\substack{x\in h\ZZ\\
    |x|\le K}}%
    \sum_{(I,j)\in S}%
    P(M_n(1+\del)>x+y,\,\At(I,j,x,h))\\
  &\le\sum_{\substack{x\in h\ZZ\\
    |x|\le K}}%
    \sum_{(I,j)\in S}\ph_{k-1}^h(x,y,\delta)P(\At(I,j,x,h))\\
  &=\sum_{\substack{x\in h\ZZ\\
    |x|\le K}}%
    \sum_{(I,j)\in S}\ph_{k-1}^h(x,y,\delta)%
    \frac{P(\At)}{P(A)}P(A(I,j,x,h)).
  \end{align*}
Note that $P(A(I,j,x,h))\ge
P(B_j(1)\in[x,x+h))\Phi(x)^{k-1}\Phi(-x-h)^{n-k}$, so that
  \[
  \frac{P(\At)}{P(A)}\le
    \left[{\frac{\Phi(x+h)}{\Phi(x)}}\right]^{k-1}
    \left[{\frac{\Phi(-x)}{\Phi(-x-h)}}\right]^{n-k}.
  \]
If we denote the right-hand side of this inequality by $g_h(x)$,
then by dominated convergence,
  \begin{align*}
  P(\Del M_n>y,\,|M_n(1)|\le K)
    &\le\sum_{\substack{x\in h\ZZ\\
    |x|\le K}}%
    \ph_{k-1}^h(x,y,\delta)g_h(x)\sum_{(I,j)\in S}P(A(I,j,x,h))\\
  &=\sum_{\substack{x\in h\ZZ\\
    |x|\le K}}%
    \ph_{k-1}^h(x,y,\delta)g_h(x)P(M_n(1)\in[x,x+h))\\
  &\to\int_{-K}^K{\ph_{k-1}(x,y,\delta)f_n(x)\,dx}.
  \end{align*}
Letting $K\to\infty$ finishes the proof. \qed

The estimate in Theorem \ref{T:cond} can be simplified even
further and we will find it convenient to use the following.

\begin{corollary}\label{C:cond}
Let $n\ge3$, $k=\flr{(n+1)/2}$, $y>0$, and $\del>0$. Then
  \begin{equation}
  P(M_n(1+\del)-M_n(1)>y)%
    \le\ph_{k-1}(x_0,y,\del)+P(M_n(1)\le x_0)\label{44a}
  \end{equation}
for all $x_0\in\RR$.
\end{corollary}

\pf We will first show that $x\mapsto\ph_{k-1}(x,y,\del)$ is
decreasing, for which it will suffice to show that $x\mapsto
p_1(x,y,\del)$ is increasing. To see this, recall that
$\ph_{k-1}(x,y,\del)=P(\sum_{j=1}^{k-1}Y_j\ge0)$. If $x\mapsto
p_1(x,y,\del)$ is increasing, then $x\mapsto
p_2(x,y,\del)=p_1(-x,-y,\del)$ is decreasing. Hence, by
\eqref{38a} and \eqref{38b}, $P(Y_j=-1)=p_1(1-p_2)$ increases with
$x$ and $P(Y_j=1)=p_2(1-p_1)$ decreases with $x$, which shows that
$x\mapsto\ph_{k-1}(x,y,\del)$ is decreasing.

With $\psi$ as in \eqref{45a}, we have $p_1=\psi/\Phi(x)$ and
  \begin{equation}
  \pa_xp_1=-\frac{\Phi'(x)}{[\Phi(x)]^2}\,\psi%
    +\frac1{\Phi(x)}\left[{%
    \Phi\left({\frac y{\sqrt{\del}}}\right)\Phi'(x)%
    +\frac1{\sqrt{\del}}\int_{-\infty}^x{%
    \Phi'\left({\frac{x+y-t}{\sqrt{\del}}}\right)\Phi'(t)\,dt}%
    }\right].\label{45a2}
  \end{equation}
Integrating by parts gives
  \[
  \psi(x,y,\del)=\Phi\left({\frac y{\sqrt{\del}}}\right)\Phi(x)%
    +\frac1{\sqrt{\del}}\int_{-\infty}^x{%
    \Phi'\left({\frac{x+y-t}{\sqrt{\del}}}\right)\Phi(t)\,dt}.
  \]
Substituting this into \eqref{45a2} gives
  \begin{align}
  \pa_xp_1%
    &=-\frac{\Phi'(x)}{[\Phi(x)]^2\sqrt{\del}}\int_{-\infty}^x{%
    \Phi'\left({\frac{x+y-t}{\sqrt{\del}}}\right)\Phi(t)\,dt}%
    +\frac1{\Phi(x)\sqrt{\del}}\int_{-\infty}^x{%
    \Phi'\left({\frac{x+y-t}{\sqrt{\del}}}\right)\Phi'(t)%
    \,dt}\notag\\
  &=\frac1{\Phi(x)\sqrt{\del}}\int_{-\infty}^x{%
    \Phi'\left({\frac{x+y-t}{\sqrt{\del}}}\right)%
    \left[{\frac{\Phi'(t)}{\Phi(t)}-\frac{\Phi'(x)}{\Phi(x)}
    }\right]\Phi(t)\,dt}.\label{45b}
  \end{align}
Note that
  \begin{align*}
  \frac d{dx}\left[{\frac{\Phi'(x)}{\Phi(x)}}\right]
    &=\frac{\Phi''(x)\Phi(x)-[\Phi'(x)]^2}{[\Phi(x)]^2}\\
  &=\frac1{[\Phi(x)]^2}\left({
    -\frac1{\sqrt{2\pi}}xe^{-x^2/2}\Phi(x)-\frac1{2\pi}e^{-x^2}
    }\right)\\
  &=-\frac{e^{-x^2/2}}{\sqrt{2\pi}[\Phi(x)]^2}%
    \left({x\Phi(x)+\frac1{\sqrt{2\pi}}e^{-x^2/2}}\right).
  \end{align*}
Clearly, $x\,\Phi(x)+\frac{1}{\sqrt{2\pi}}\,e^{-x^2/2}\ge 0$ for
$x\ge 0$. If $x<0$, then by \eqref{25b},
  \begin{align*}
  x\Phi(x)+\frac1{\sqrt{2\pi}}e^{-x^2/2}
    &=x\Phi(-|x|)+\frac1{\sqrt{2\pi}}e^{-x^2/2}\\
  &\ge x\frac1{\sqrt{2\pi}}|x|^{-1}e^{-x^2/2}
    +\frac1{\sqrt{2\pi}}e^{-x^2/2}\\
  &=0.
  \end{align*}
Thus, $x\mapsto\Phi'(x)/\Phi(x)$ is decreasing, so by \eqref{45b},
$\pa_xp_1\ge 0$.

Hence, $x\mapsto\ph_{k-1}(x,y,\del)$ is decreasing, and using
Theorem \ref{T:cond},
  \begin{align*}
  P(M_n(1+\del)-M_n(1)>y)
    &\le\int_{-\infty}^\infty{\ph_{k-1}(x,y,\del)f_n(x)\,dx}\\
  &\le\int_{-\infty}^{x_0}{\ph_{k-1}(x,y,\del)f_n(x)\,dx}
    +\ph_{k-1}(x_0,y,\del)\int_{x_0}^\infty{f_n(x)\,dx}\\
  &\le\int_{-\infty}^{x_0}{f_n(x)\,dx}
    +\ph_{k-1}(x_0,y,\del)\int_{-\infty}^\infty{f_n(x)\,dx}\\
  &=P(M_n(1)\le x_0)+\ph_{k-1}(x_0,y,\del),
  \end{align*}
where $x_0\in\RR$ is arbitrary. \qed

Recall that our only remaining goal is to establish the inequality
\eqref{21a} for the small and medium jump regimes. In applying
Corollary \ref{C:cond} to this task, we must set $y=\ee/\sqrt{n}$.
Our choice for $x_0$, however, is less clear. On the one hand, we
want $x_0$ to be large so that the first term on the right-hand of
\eqref{44a} is small. On the other hand, we need $x_0$ to be
sufficiently far into the negative real line so that the second
term is small. The value of $x_0$ that will strike a balance for
us is given in the following lemma.

\begin{lemma}\label{L:balance}
Let $\ee>0$, $\del>0$, and $n\in\NN$. Define
$x_0=-\ee/(\del^{1/4}\sqrt{n})$. Then for all $p>2$,
  \[
  P(M_n(1)\le x_0)\le C_p(\ee^{-1}\del^{1/4})^p,
  \]
where $C_p$ is a finite constant depending only $p$.
\end{lemma}

\pf This follows immediately from \eqref{24a}. \qed

In light of this lemma and Corollary \ref{C:cond}, we will
establish inequality \eqref{21a} once we verify that
  \begin{equation}
  \ph_{k-1}\left({
    -\frac\ee{\del^{1/4}\sqrt{n}},\frac\ee{\sqrt{n}},\del
    }\right)\le C_p(\ee^{-1}\del^{1/6})^p\label{49a}
  \end{equation}
for all values of $\ee$, $\del$, and $n$ in the small and medium
jump regimes.

\section{Estimates for a Random Walk}

In this section, we wish to find useful estimates for
$\ph_k(x,y,\del)=P(S_k\ge 0)$. The process $\{S_n\}_{n=1}^\infty$
is, of course, a biased random walk which, in the cases we are
interested in, has a negative drift. Let us recall the definition
of $S_n$. In this section, we will temporarily abandon the tilde
notation for the sake of simplicity.

We take as given a sequence of $\{-1,0,1\}$-valued random
variables with $p_1=P(Y_j=-1)$ and $p_2=P(Y_j=1)$. We define
$\ee=p_1+p_2$ and $\mu=p_1-p_2$, so that $P(Y_j=0)=1-\ee$. We then
define $S_n=\sum_{j=1}^n Y_j$.

As mentioned, we will be interested in the case where $\mu>0$, so
that the walk has a negative drift. Besides this, however, we will
also be interested in the case where $\ee$ is small. That is,
besides the negative drift, our walk will have the property that,
for most time steps, it does not move. Our first estimate is a
straightforward application of Chebyshev's inequality. It is a
fairly simple result and serves as our starting point, but it will
not be sufficient by itself. Note, in particular, that it does not
make any noteworthy use of the fact that $\ee$ is small.

\begin{lemma}\label{L:cheby}
If $\ee>0$ and $\mu>0$, then for all $p>1$, there exists $C_p$,
depending only on $p$, such that
  \begin{equation}
  P(S_n\ge 0)\le C_p\frac\ee{n^p\mu^{2p}}\label{53a}
  \end{equation}
for all $n$.
\end{lemma}

\pf Since $EY_j=-\mu$, Chebyshev's inequality gives
  \[
  P(S_n\ge 0)=P(S_n+n\mu\ge n\mu)
    \le\frac{E|S_n+n\mu|^{2p}}{n^{2p}\mu^{2p}}.
  \]
By Burkholder's and Jensen's inequalities,
  \[
  E|S_n+n\mu|^{2p}=E\biggl|{\sum_{j=1}^n(Y_j+\mu)}\biggr|^{2p}
    \le\Ct_p E\biggl|{\sum_{j=1}^n|Y_j+\mu|^2}\biggr|^p
    \le\Ct_p n^p E|Y_1+\mu|^{2p}.
  \]
Also,
  \begin{align*}
  E|Y_1+\mu|^{2p}
    &=p_1(1-\mu)^{2p}+(1-\ee)\mu^{2p}+p_2(1+\mu)^{2p}\\
  &\le 2^{2p}(p_1+p_2)+\mu^{2p}\\
  &\le(2^{2p}+1)\ee
  \end{align*}
since $\mu\le\ee$. Thus, \eqref{53a} holds with
$C_p=\Ct_p(2^{2p}+1)$. \qed

As it stands, \eqref{53a} will not suit our needs. We will find it
necessary for the numerator on the right-hand side of \eqref{53a}
to contain $\ee^p$ rather than $\ee$. To accomplish this, we must
appeal to the fact that, for the most part, this random walk does
not move. To this end, we begin with two lemmas.

\begin{lemma}\label{L:gauss}
For $n\in\NN$, $k\in\{0,\ldots,n\}$, $p\in(0,1)$, and $x\in\RR$,
let $f(n,k,p)=\binom{n}{k}p^k q^{n-k}$, where $q=1-p$, and let
$g(n,x,p)=(2\pi npq)^{-1/2}\exp\{-(x-np)^2/2npq\}$. Then
  \[
  \sup_{n\in\NN}\left({\sup_{k\in\{0,\ldots,n\}}
    \frac{f(n,k,p)}{g(n,k,p)}}\right)<\infty
  \]
if and only if $p=1/2$. However, there exists a universal
constant $C$, independent of $p$, such that %
$f(n,k,p)/g(n,k,p)\le C$ for all $n\in\NN$ and all
$k\in\{0,\ldots,\flr{np}\}$, provided $p\le 1/2$.
\end{lemma}

\pf It will first be shown that there exists a universal constant
$C$ such that
  \begin{description}
  \item{(i)} if $p\le 1/2$, then $f(n,0,p)/g(n,0,p)\le C$, and%
  \item{(ii)} if $p\le 1/2$ and $\flr{np}\ge 1$, then %
    $f(n,1,p)/g(n,1,p)\le C$.
  \end{description}
We will start by showing that if $\al>0$, then there exists a
constant $C_{\al}$, depending only on $\al$, such that for all
$p\le 1/2$,
  \begin{equation}
  (np)^{\al}(qe^{p/2q})^n\le C_{\al}.\label{54a}
  \end{equation}
To prove this, first consider $2/5\le p\le 1/2$. In this case,
$qe^{p/2q}\le\frac35e^{1/2}<1$. Thus,
  \[
  (np)^{\al}(qe^{p/2q})^n\le\sup_n\left[{n^{\al}
    \left({\frac35e^{1/2}}\right)^n}\right]<\infty.
  \]
Next, consider $0<p<2/5$. Since %
$\frac d{dq}[\log(q^{5/6}e^{p/2q})]=(5q-3)/6q^2>0$ for $q>3/5$, it
follows that in this case, $q^{5/6}e^{p/2q}\le 1$. Hence,
  \[
  (np)^{\al}(qe^{p/2q})^n\le(np)^{\al}q^{n/6}
    =(n^{\al}q^{n/6})p^{\al}.
  \]
Elementary calculus shows that $x\mapsto x^{\al}q^{x/6}$ attains
its maximum on $[0,\infty)$ at $x=-6\al/\log q$. Thus,
  \[
  n^{\al}q^{n/6}p^{\al}\le\left({\frac{6\alpha}e}\right)^{\alpha}
    \left({\frac{1-q}{|\log q|}}\right)^{\alpha}.
  \]
Since $(q-1)/\log q\to 1$ as $q\to 1$, this proves \eqref{54a}.
Thus, if $p\le 1/2$, then
  \[
  \frac{f(n,0,p)}{g(n,0,p)}=\sqrt{2\pi npq}\,q^n e^{np/2q}
    =\sqrt{2\pi q}(np)^{1/2}(qe^{p/2q})^n
    \le\sqrt{2\pi}C_{1/2},
  \]
and if $p\le 1/2$ and $np\ge 1$, then
  \begin{align*}
  \frac{f(n,1,p)}{g(n,1,p)}&=\sqrt{2\pi npq}\,npq^{n-1}
    \exp\left\{{\frac{np}{2q}-\frac1q+\frac1{2npq}}\right\}\\
  &\le\sqrt{2\pi q}\,q^{n-1}(np)^{3/2}e^{np/2q}\\
  &=\sqrt{\frac{2\pi}q}(np)^{3/2}(qe^{p/2q})^n\\
  &\le\sqrt{4\pi}C_{3/2},
  \end{align*}
which verifies (i) and (ii).

Now, for $k\in\{1,\ldots,n-1\}$, Stirling's formula implies that
$f(n,k,p)$ is bounded above and below by universal, positive
constant multiples of
  \[
  \frac{n^{n+\frac12}}{(n-k)^{n-k+\frac12}k^{k+\frac12}}\,
    p^kq^{n-k}.
  \]
Let us define
  \begin{align*}
  F(k)=F(n,k,p)&=\log\left({
    \frac{n^{n+\frac12}}{(n-k)^{n-k+\frac12}k^{k+\frac12}}\,
    p^kq^{n-k}}\right)-\log(\sqrt{2\pi}\,g(n,k,p))\\
  &=(n+1)\log n-(n-k+\tfrac12)\log(n-k)-(k+\tfrac12)\log k\\
  &\quad+(k+\tfrac12)\log p
    +(n-k+\tfrac12)\log q+(k-np)^2/2npq,
  \end{align*}
so that there are universal, positive constants $C_1$ and $C_2$
such that
  \begin{equation}
  \log C_1+F(k)\le\log\left[{\frac{f(n,k,p)}{g(n,k,p)}}\right]
    \le\log C_2+F(k)\label{56a}
  \end{equation}
for all $k\in\{1,\ldots,n-1\}$. Note that $F(k)$ is well-defined
for all real $k\in(0,n)$.

We can directly compute that
  \[
  F(n/2)=\frac12\log(4pq)+\frac n2(G(p)+G(1-p)),
  \]
where $G(p)=\log 2+\log p+1/(4p)-1/2$. Now, $G'(p)=1/p-1/(4p^2)$,
which gives
  \[
  G'(p)-G'(1-p)=\left({\frac{q-p}{pq}}\right)
    \left({1-\frac{1}{4pq}}\right).
  \]
Since $1-1/(4pq)<0$ for all $p\ne 1/2$, the function %
$p\mapsto G(p)+G(1-p)$ is strictly decreasing on $(0,1/2)$ and
strictly increasing on $(1/2,0)$. Since $G(1/2)=0$, we have that
$G(p)+G(1-p)>0$ for all $p\ne 1/2$. Thus, if $p\ne 1/2$, then
$F(n/2)\to\infty$ as $n\to\infty$. It now follows from \eqref{56a}
that
  \[
  \sup_{n\in\NN}\left({\sup_{k\in\{0,\ldots,n\}}
    \frac{f(n,k,p)}{g(n,k,p)}}\right)=\infty
  \]
whenever $p\ne 1/2$.

Now suppose $p\le 1/2$ and let $k\in[2,np]$. We can compute that
for all $x\in(0,n)$,
  \begin{align*}
  F'(x)&=\log(n-x)+\frac1{2(n-x)}-\log x-\frac1{2x}\log\frac pq
    +\frac{x}{npq}-\frac1q\\
  F''(x)&=-\frac1{n-x}+\frac1{2(n-x)^2}-\frac1x
    +\frac1{2x^2}+\frac1{npq}\\
  F'''(x)&=-\frac1{(n-x)^2}+\frac1{(n-x)^3}+\frac1{x^2}
    -\frac1{x^3}\\
  F^{(4)}(x)&=\frac{3-2(n-x)}{(n-x)^4}+\frac{3-2x}{x^4}.
  \end{align*}
It is easily verified that $F(np)=0$ and $F'(np)=(p-q)/2npq$, so
that we may write
  \begin{align*}
  F(k)&=-\int_k^{np}{F'(t)\,dt}\\
  &=-\int_k^{np}{\left({
    \frac{p-q}{2npq}-\int_t^{np}{F''(s)\,ds}}\right)\,dt}\\
  &\le\frac{q-p}{2q}+\int_k^{np}{\int_k^s{F''(s)\,dt}\,ds}.
  \end{align*}
Since $F^{(4)}\le 0$ on $[2,n-2]$ and $F'''(n/2)=0$, it follows
that $F'''\ge 0$ on $[2,n/2]$, which implies $F''$ is increasing
on $[2,n/2]$. Since $F''(np)=(p^2+q^2)/2n^2p^2q^2$, we have
  \[
  F(k)\le\frac12+\frac{p^2+q^2}{2n^2p^2q^2}
    \int_k^{np}{(s-k)\,ds}
    \le\frac12+\frac{p^2+q^2}{2n^2p^2q^2}n^2p^2
    \le\frac32
  \]
for all $p\le 1/2$.

It now follows from \eqref{56a} and (i), (ii) that there is a
universal constant $C$, independent of $p$, such that
$f(n,k,p)/g(n,k,p)\le C$ for
all $n\in\NN$ and all $k\in\{0,\ldots,\flr{np}\}$, provided %
$p\le 1/2$. Also, if $p=1/2$, symmetry gives the same bound for
$k\in\{\flr{n/2}+1,\ldots,n\}$, and it follows that
  \[
  \sup_{n\in\NN}\left({\sup_{k\in\{0,\ldots,n\}}
    \frac{f(n,k,p)}{g(n,k,p)}}\right)<\infty,
  \]
which completes the proof. \qed

\begin{lemma}\label{L:recip}
Let $0<\ee<1/2$ and suppose that $\{\xi_j\}_{j=1}^\infty$ are
i.i.d. $\{0,1\}$-valued random variables with $P(\xi_1=1)=\ee$.
Let $T_n=\sum_{j=1}^n\xi_j$. Then for each $p>1$, there exists a
finite constant $C_p$, depending only on $p$, such that
  \[
  E[T_n^{-p}1_{\{T_n>0\}}]\le C_p\frac1{(\ee n)^p}
  \]
for all $n\in\NN$.
\end{lemma}

\pf Observe that
  \begin{align*}
  E[T_n^{-p}1_{\{T_n>0\}}]&=E[T_n^{-p}1_{\{1\le T_n\le\ee n/2\}}]
    +E[T_n^{-p}1_{\{T_n>\ee n/2\}}]\\
  &\le P\left({T_n\le\frac{\ee n}2}\right)
    +\left({\frac{\ee n}2}\right)^{-p}.
  \end{align*}
Hence, it will suffice to show that
  \[
  P\left({T_n\le\frac{\ee n}2}\right)
    \le C_p\frac1{(\ee n)^p}.
  \]
To see this, let $f$ and $g$ be as in Lemma \ref{L:gauss} with
$p=\ee$, so that there exists a universal, finite constant $C$,
independent of $\ee$, such that $f(n,k,\ee)\le Cg(n,k,\ee)$ for
all $n\in\NN$ and all $k\in\{0,\ldots,\flr{\ee n}\}$. Let
$m=\flr{\ee n/2}$, so that
  \[
  P\left({T_n\le\frac{\ee n}2}\right)=P(T_n\le m)
    =\sum_{k=0}^m P(T_n=k)\le C\sum_{k=0}^m g(n,k,\ee).
  \]
If $\ee n\le 4$, then $P(T_n\le m)\le 1\le 4^p/(\ee n)^p$, so that
we may assume without loss of generality that $\ee n>4$. Note that
$x\mapsto g(n,x,\ee)$ is increasing on $[0,\ee n]$ and $\ee n>4$
implies $m+1\le(\ee n/2)+1<3\ee n/4$. Thus,
  \begin{align*}
  P(T_n\le m)&\le C\int_0^{m+1}{g(n,x,\ee)\,dx}\\
  &\le C\int_{-\infty}^{3\ee n/4}{g(n,x,\ee)\,dx}\\
  &=\frac C{\sqrt{2\pi t}}
    \int_{-\infty}^{3\ee n/4}{e^{-(x-\ee n)^2/2t}\,dx},
  \end{align*}
where $t=n\ee(1-\ee)$. By a change of variables,
  \[
  P(T_n\le m)\le C\Phi\left({-\frac{\ee n}{4\sqrt{t}}}\right)
    \le C\Phi\left({-\frac{\sqrt{\ee n}}4}\right).
  \]
By \eqref{25b},
  \[
  P(T_n\le m)\le\frac C{\sqrt{2\pi}}
    \cdot\frac4{\sqrt{\ee n}}\,e^{-\ee n/32}
    \le C\sqrt{\frac2\pi}\,e^{-\ee n/32}.
  \]
Since there exists $K_p<\infty$ such that $x^p e^{-x/32}\le K_p$
for all $x\in[0,\infty)$, we have
  \[
  P(T_n\le m)\le C\sqrt{\frac2\pi}\,K_p\frac1{(\ee n)^p},
  \]
which finishes the proof. \qed

With these lemmas in place, we may now make the needed improvement
to Lemma \ref{L:cheby}.

\begin{lemma}\label{L:chebyplus}
If $0<\ee<1/2$ and $\mu>0$, then for all $p>1$, there exists
$C_p$, depending only on $p$, such that
  \begin{equation}
  P(S_n\ge 0)\le C_p\frac{\ee^p}{n^p\mu^{2p}}\label{61a}
  \end{equation}
for all $n$.
\end{lemma}

\pf Let $\{\Yt_j\}_{j=1}^\infty$ be a sequence of i.i.d.
$\{-1,1\}$-valued random variables with $P(\Yt_1=-1)=p_1/\ee$. Let
$\{\xi_j\}_{j=1}^\infty$ be a sequence of i.i.d. $\{0,1\}$-valued
random variables, independent of $\{\Yt_j\}_{j=1}^\infty$, with
$P(\xi_1=1)=\ee$. Then $\{\Yt_j\xi_j\}_{j=1}^\infty$ is an i.i.d.
sequence of random variables which has the same law as
$\{Y_j\}_{j=1}^\infty$.

Let $\St_n=\sum_{j=1}^n \Yt_j$ and note that by Lemma
\ref{L:cheby},
  \begin{equation}
  P(\St_n\ge 0)\le\Ct_p\frac1{n^p(\mu/\ee)^{2p}}
    =\Ct_p\frac{\ee^{2p}}{n^p\mu^{2p}}.\label{61b}
  \end{equation}
Define $\xi^{(n)}=(\xi_1,\ldots,\xi_n)$, so that
  \begin{align*}
  P(S_n\ge 0)&=P\biggl({\sum_{j=1}^n\Yt_j\xi_j\ge 0}\biggr)\\
  &=\sum_{k=0}^n\sum_{\substack{\al\in\{0,1\}^n\\
    |\al|=k}}
    P\biggl({\sum_{j=1}^n\Yt_j\xi_j\ge 0,\,\xi^{(n)}=\al}\biggr)\\
  &=\sum_{k=0}^n\sum_{\substack{\al\in\{0,1\}^n\\
    |\al|=k}}
    P\biggl({\sum_{\{j:\al_j=1\}}
    \Yt_j\ge 0,\,\xi^{(n)}=\al}\biggr),
  \end{align*}
where $|\al|=\al_1+\cdots+\al_n$. If $T_n=\sum_{j=1}^n\xi_j$, then
by symmetry and independence,
  \[
  P(S_n\ge 0)=\sum_{k=0}^n\sum_{\substack{\al\in\{0,1\}^n\\
    |\al|=k}}
    P\biggl({\sum_{j=1}^k\Yt_j\ge 0}\biggr)P(\xi^{(n)}=\al)
    =\sum_{k=0}^nP(\St_k\ge 0)P(T_n=k).
  \]
Using \eqref{61b} and Lemma \ref{L:recip},
  \begin{align*}
  P(S_n\ge 0)&\le P(T_n=0)+\Ct_p\frac{\ee^{2p}}{\mu^{2p}}\,
    \sum_{k=1}^n k^{-p}P(T_n=k)\\
  &=(1-\ee)^n+\Ct_p\frac{\ee^{2p}}{\mu^{2p}}\,
    E[T_n^{-p}1_{\{T_n>0\}}]\\
  &\le(1-\ee)^n+\Ct'_p\frac{\ee^{2p}}{\mu^{2p}}\,
    \frac1{(\ee n)^p},
  \end{align*}
Note that $1-\ee\le e^{-\ee}$, so that
  \[
  (1-\ee)^n\le e^{-\ee n}\le\Ct''_p\frac1{(\ee n)^p}
    =\Ct''_p\frac{\ee^p}{n^p\ee^{2p}}
    \le\Ct''_p\frac{\ee^p}{n^p\mu^{2p}},
  \]
which gives \eqref{61a} with $C_p=\Ct''_p+\Ct'_p$. \qed

\section{The Small Jump Regime}\label{small}

Let us now put the pieces together and establish \eqref{21a} for
the small jump regime. Recall from Section \ref{cond} that it will
suffice to establish \eqref{49a}. Using the notation of
\eqref{37a}-\eqref{38d}, Lemma \ref{L:chebyplus} will give us
that, for $p>1$,
  \begin{equation}
  \ph_{k-1}(x,y,\del)
    \le C_p\frac{\eet^p}{(k-1)^p\mut^{2p}},\label{63a}
  \end{equation}
provided $\eet=\eet(x,y,\del)<1/2$ and $\mut=\mut(x,y,\del)>0$. We
will be applying this with $x=-\ee/(\del^{1/4}\sqrt{n})$ and
$y=\ee/\sqrt{n}$, but recall that in the small jump regime, we can
write $\ee/\sqrt{n}=\del^{1/2+\al}$ for some $\al>0$. As such, the
following lemma will help us check the provisions of \eqref{63a}.

\begin{lemma}\label{L:taylor}
For each $\Del>0$, there exists $\del_0>0$ such that
  \begin{description}
  \item{(i)} $\mut(-\del^{1/4+\al},\del^{1/2+\al},\del)
    \ge\frac1{\sqrt{2\pi}}\,\del^{1/2+\al}$, and
  \item{(ii)} $\eet(-\del^{1/4+\al},\del^{1/2+\al},\del)
    \le 1000\,\del^{1/2}<\frac12$
  \end{description}
for all $\al\ge\Del$ and all $0<\del\le\del_0$.
\end{lemma}

\pf For fixed $\del>0$, let $\psi(x,y)=\psi(x,y,\del)$ be given by
\eqref{45a}. We wish to show that
  \begin{equation}
  \psi(x,y)=\frac12-\frac1{2\pi}\tan^{-1}\sqrt{\del}
    +\frac x{\sqrt{2\pi}}+\frac y{2\sqrt{2\pi}}
    +\frac{\sqrt{\del}}{4\pi}(x+y)^2
    -\frac{y^2}{4\pi\sqrt{\del}}+\Rt(x,y),\label{65c}
  \end{equation}
where
  \begin{equation}
  |\Rt(x,y)|\le(|x|+|y|)^3+\frac{|x||y|^2}{\sqrt{\del}}(|x|+|y|)
    +\frac{|y|^4}{\del^{3/2}}+\del^{3/2}(x+y)^2
    +\del(|x|+|y|)\label{65d}
  \end{equation}
for all $x,y\in\RR$.

We will first show that for $i\ge 0$ and $j\ge 1$,
  \begin{align}
  \pa_x^i\psi&=\int_{-\infty}^x{
    \Phi\left({\frac{x+y-t}{\sqrt{\del}}}\right)
    \Phi^{(i+1)}(t)\,dt},\label{65a}\\
  \pa_x^i\pa_y^j\psi&=-\left({\frac1{\sqrt{\del}}}\right)^{j-1}
    \Phi^{(j-1)}\left({\frac y{\sqrt{\del}}}\right)
    \Phi^{(i+1)}(x)+\pa_x^{i+1}\pa_y^{j-1}\psi.\label{65b}
  \end{align}
For $i=0$, \eqref{65a} is just the definition of $\psi$. If
\eqref{65a} is true for some $i\ge 0$, then using integration by
parts gives
  \begin{align*}
  \pa_x^{i+1}\psi &=\pa_x\left[{\int_{-\infty}^x{
    \Phi\left({\frac{x+y-t}{\sqrt{\del}}}\right)
    \Phi^{(i+1)}(t)\,dt}}\right]\\
  &=\Phi\left({\frac y{\sqrt{\del}}}\right)\Phi^{(i+1)}(x)
    +\frac1{\sqrt{\del}}\int_{-\infty}^x{
    \Phi'\left({\frac{x+y-t}{\sqrt{\del}}}\right)
    \Phi^{(i+1)}(t)\,dt}\\
  &=\int_{-\infty}^x{
    \Phi\left({\frac{x+y-t}{\sqrt{\del}}}\right)
    \Phi^{(i+2)}(t)\,dt},
  \end{align*}
so by induction, \eqref{65a} holds for all $i\ge 0$. For
\eqref{65b}, first consider $j=1$. Then
  \begin{align*}
  \pa_x^i\pa_y\psi &=\pa_y\left[{
    \int_{-\infty}^x{\Phi\left({\frac{x+y-t}{\sqrt{\del}}}\right)
    \Phi^{(i+1)}(t)\,dt}}\right]\\
  &=\int_{-\infty}^x{\pa_y\left[{
    \Phi\left({\frac{x+y-t}{\sqrt{\del}}}\right)
    }\right]\Phi^{(i+1)}(t)\,dt}\\
  &=\int_{-\infty}^x{\pa_x\left[{
    \Phi\left({\frac{x+y-t}{\sqrt{\del}}}\right)
    }\right]\Phi^{(i+1)}(t)\,dt}\\
  &=\pa_x\left[{
    \int_{-\infty}^x{\Phi\left({\frac{x+y-t}{\sqrt{\del}}}\right)
    \Phi^{(i+1)}(t)\,dt}}\right]
    -\Phi\left({\frac y{\sqrt{\del}}}\right)\Phi^{(i+1)}(x)\\
  &=-\Phi\left({\frac y{\sqrt{\del}}}\right)\Phi^{(i+1)}(x)
    +\pa_x^{i+1}\psi,
  \end{align*}
and \eqref{65b} holds for all $i\ge 0$ when $j=1$. Now suppose
\eqref{65b} holds for some $j\ge 1$ and all $i\ge 0$. Then
  \begin{align*}
  \pa_x^i\pa_y^{j+1}\psi &=\pa_y\left[{
    -\left({\frac1{\sqrt{\del}}}\right)^{j-1}
    \Phi^{(j-1)}\left({\frac y{\sqrt{\del}}}\right)
    \Phi^{(i+1)}(x)+\pa_x^{i+1}\pa_y^{j-1}\psi}\right]\\
  &=-\left({\frac1{\sqrt{\del}}}\right)^j
    \Phi^{(j)}\left({\frac y{\sqrt{\del}}}\right)
    \Phi^{(i+1)}(x)+\pa_x^{i+1}\pa_y^j\psi.
  \end{align*}
By induction, \eqref{65b} holds for all $i\ge 0$ and $j\ge 1$.

By Taylor's Theorem we have that
  \begin{equation}\begin{split}
  \psi(x,y)&=\psi(0,0)+x\psi_x(0,0)+y\psi_y(0,0)\\
  &\quad+\frac1{2!}[x^2\psi_{xx}(0,0)+2xy\psi_{xy}(0,0)
    +y^2\psi_{yy}(0,0)]+R^{(1)}(x,y),\label{66a}
  \end{split}\end{equation}
where
  \[
  R^{(1)}(x,y)=\frac1{3!}[x^3\psi_{xxx}(\bar{x},\bar{y})
    +3x^2y\psi_{xxy}(\bar{x},\bar{y})
    +3xy^2\psi_{xyy}(\bar{x},\bar{y})
    +y^3\psi_{yyy}(\bar{x},\bar{y})]
  \]
and $(\bar{x},\bar{y})=(\theta x,\theta y)$ for some
$\th\in(0,1)$. Using \eqref{65a}, \eqref{65b}, and direct
integration, we can verify that \eqref{66a} becomes
  \begin{equation*}\begin{split}
  \psi(x,y)&=\frac12-\frac1{2\pi}\tan^{-1}\sqrt{\del}
    +\frac x{2\sqrt{2\pi}}\left({1+\frac1{\sqrt{1+\del}}}\right)
    +\frac y{2\sqrt{2\pi}\sqrt{1+\del}}\\
  &\quad+\frac{(x+y)^2\sqrt{\del}}{4\pi(1+\del)}
    -\frac{y^2}{4\pi\sqrt{\del}}+R^{(1)}(x,y).
  \end{split}\end{equation*}
Now,
  \begin{align*}
  \frac x{2\sqrt{2\pi}}\left({1+\frac1{\sqrt{1+\del}}}\right)
    &=\frac x{\sqrt{2\pi}}
    +\frac x{2\sqrt{2\pi}}\left({\frac1{\sqrt{1+\del}}-1}\right)\\
  \frac y{2\sqrt{2\pi}\sqrt{1+\del}}&=\frac y{2\sqrt{2\pi}}
    +\frac{y}{2\sqrt{2\pi}}
    \left({\frac1{\sqrt{1+\del}}-1}\right)\\
  \frac{(x+y)^2\sqrt{\del}}{4\pi(1+\del)}
    &=\frac{\sqrt{\del}}{4\pi}(x+y)^2
    +\frac{\sqrt{\del}}{4\pi}(x+y)^2
    \left({\frac1{1+\del}-1}\right).
  \end{align*}
Thus, if
  \[
  R^{(2)}(x,y)=\frac{x+y}{2\sqrt{2\pi}}
    \left({\frac1{\sqrt{1+\del}}-1}\right)
    -\frac{\del^{3/2}(x+y)^2}{4\pi(1+\del)},
  \]
then \eqref{65c} holds with $\Rt=R^{(1)}+R^{(2)}$.

Since $|(1+\del)^{-1/2}-1|<\del$, we have
$|R^{(2)}(x,y)|\le\del(|x|+|y|)+\del^{3/2}(x+y)^2$. To estimate
$R^{(1)}$, we must estimate the third partial derivatives of
$\psi$. Using \eqref{65a}, we have
  \[
  |\psi_{xxx}(x,y)|=\left|{\int_{-\infty}^x{
    \Phi\left({\frac{x+y-t}{\sqrt{\del}}}\right)
    \Phi^{(4)}(t)\,dt}}\right|
    \le\int_{-\infty}^\infty{|\Phi^{(4)}(t)|\,dt}.
  \]
Since $\Phi^{(4)}(t)=(3t-t^3)\Phi'(t)$, we have
  \[
  |\psi_{xxx}(x,y)|\le 2\int_0^\infty{(3t+t^3)\Phi'(t)\,dt}
    =\frac{10}{\sqrt{2\pi}}.
  \]
Similarly, by \eqref{65b},
  \[
  |\psi_{xxy}(x,y)|
    =\left|{-\Phi\left({\frac{y}{\sqrt{\del}}}\right)\Phi'''(x)
    +\psi_{xxx}(x,y)}\right|
    \le|\Phi'''(x)|+\frac{10}{\sqrt{2\pi}}.
  \]
Since $|\Phi'''(x)|\le 2(2\pi)^{-1/2}$ for all $x\in\RR$, we have
that
  \[
  |\psi_{xxy}(x,y)|\le\frac{12}{\sqrt{2\pi}}.
  \]
Likewise, the formulas
  \begin{align*}
  \psi_{xyy}
    &=-\frac1{\sqrt{\del}}\Phi'\left({\frac y{\sqrt{\del}}}\right)
    \Phi''(x)+\psi_{xxy}\\
  &=\frac x{\sqrt{\del}}\Phi'\left({\frac y{\sqrt{\del}}}\right)
    \Phi'(x)+\psi_{xxy}
  \end{align*}
and
  \begin{align*}
  \psi_{yyy}
    &=-\frac1{\del}\Phi''\left({\frac y{\sqrt{\del}}}\right)
    \Phi'(x)+\psi_{xyy}\\
  &=\frac y{\del^{3/2}}\Phi'\left({\frac y{\sqrt{\del}}}\right)
    \Phi'(x)+\psi_{xyy}
  \end{align*}
can be used to verify that
  \begin{align*}
  |\psi_{xyy}(x,y)|&\le(|x|\del^{-1/2}+12\sqrt{2\pi})/(2\pi)\\
  |\psi_{yyy}(x,y)|&\le(|y|\del^{-3/2}+|x|\del^{-1/2}
    +12\sqrt{2\pi})/(2\pi).
  \end{align*}
Piecing this together, we have
  \begin{align*}
  |R^{(1)}(x,y)|&\le\frac1{3!}\left[{\frac{10|x|^3}{\sqrt{2\pi}}
    +\frac{36|x|^2|y|}{\sqrt{2\pi}}
    +3|x||y|^2\left({\frac{|x|}{2\pi\sqrt{\del}}
    +\frac{12}{\sqrt{2\pi}}}\right)}\right.\\
  &\quad\left.{+|y|^3\left({\frac{|y|}{2\pi\del^{3/2}}
    +\frac{|x|}{2\pi\sqrt{\del}}+\frac{12}{\sqrt{2\pi}}}\right)
    }\right]\\
  &\le\frac1{3!}\left[{\frac{12}{\sqrt{2\pi}}(|x|+|y|)^3
    +\frac{3|x||y|^2}{2\pi\sqrt{\del}}(|x|+|y|)
    +\frac{|y|^4}{2\pi\del^{3/2}}}\right]\\
  &\le(|x|+|y|)^3+\frac{|x||y|^2}{\sqrt{\del}}(|x|+|y|)
    +\frac{|y|^4}{\del^{3/2}}.
  \end{align*}
Combined with the estimate for $R^{(2)}$, this verifies
\eqref{65d}.

Now, observe that $p_1(x,y,\del)=\psi(x,y)/\Phi(x)$. Write
$\Phi(x)=\frac12+\frac x{\sqrt{2\pi}}+r_1(x)$, where
$r_1(x)=\frac12x^2\Phi''(\bar{x})$ and $\bar{x}=\th x$ for some
$\th\in(0,1)$. Note that $|r_1(x)|\le\frac1{2\sqrt{2\pi}}|x|^3$.
For $x\ne-\sqrt{\pi/2}$, write %
$\Phi(x)^{-1}=(\frac12+\frac x{\sqrt{2\pi}})^{-1}+r_2(x)$, where
$r_2(x)=-r_1(x)\Phi(x)^{-1}%
(\frac12+\frac x{\sqrt{2\pi}})^{-1}$. Similarly, we may write
$\Phi(x)^{-1}=2+r_3(x)$, where
  \begin{align*}
  r_3(x)&=r_2(x)+\frac1{\frac12+\frac x{\sqrt{2\pi}}}-2\\
  &=r_2(x)-\frac{4x}{\sqrt{2\pi}+2x}.
  \end{align*}
Let us now assume $|x|\le 1$. Then $x\ne-\sqrt{\pi/2}$ and the
above applies. Note that
  \[
  |r_2(x)|\le\frac{|r_1(x)|}{\Phi(-1)\left({
    \frac12-\frac1{\sqrt{2\pi}}}\right)}
  \]
Since $\Phi(-1)\ge\frac12-\frac1{\sqrt{2\pi}}\ge\frac1{10}$, we
have $|r_2(x)|\le 100|r_1(x)|\le\frac{50}{\sqrt{2\pi}}|x|^3$.
Also,
  \[
  |r_3(x)|\le|r_2(x)|+\left({\frac4{\sqrt{2\pi}-2}}\right)|x|
    \le\frac{50}{\sqrt{2\pi}}|x|^3+\frac{20}{\sqrt{2\pi}}|x|.
  \]
Since $|x|\le 1$, this gives
$|r_3(x)|\le\frac{70}{\sqrt{2\pi}}|x|$. Applying \eqref{65c}
yields
  \begin{align}
  p_1(x,y,\del)&=\psi(x,y)\Phi(x)^{-1}\notag\\
  &=\left({\frac12+\frac x{\sqrt{2\pi}}}\right)
    \biggl({\left({\frac12+\frac x{\sqrt{2\pi}}}\right)^{-1}+r_2(x)
    }\biggr)\notag\\
  &\quad+\left({-\frac1{2\pi}\tan^{-1}\sqrt{\del}
    +\frac y{2\sqrt{2\pi}}+\frac{\sqrt{\del}}{4\pi}(x+y)^2
    -\frac{y^2}{4\pi\sqrt{\del}}}\right)(2+r_3(x))\notag\\
  &\quad+\Rt(x,y)\Phi(x)^{-1}\notag\\
  &=1-\frac1{\pi}\tan^{-1}\sqrt{\del}+\frac y{\sqrt{2\pi}}
    +\frac{\sqrt{\del}}{2\pi}(x+y)^2-\frac{y^2}{2\pi\sqrt{\del}}
    +R_\del(x,y),\label{70a}
  \end{align}
where
  \begin{align*}
  |R_\del(x,y)|&\le|r_2(x)|
    +\left({\frac{\tan^{-1}\sqrt{\del}}{2\pi}
    +\frac{|y|}{2\sqrt{2\pi}}+\frac{\sqrt{\del}}{4\pi}(x+y)^2
    +\frac{y^2}{4\pi\sqrt{\del}}}\right)|r_3(x)|
    +\frac{|\Rt(x,y)|}{\Phi(-1)}\\
  &\le\frac{50}{\sqrt{2\pi}}|x|^3
    +\left({\frac{\sqrt{\del}}{2\pi}+\frac{|y|}{2\sqrt{2\pi}}
    +\frac{\sqrt{\del}}{4\pi}(x+y)^2
    +\frac{y^2}{4\pi\sqrt{\del}}}\right)\frac{70}{\sqrt{2\pi}}|x|
    +10|\Rt(x,y)|.
  \end{align*}
Hence,
  \begin{equation*}\begin{split}
  |R_\del(x,y)|&\le\frac{50}{\sqrt{2\pi}}|x|^3
    +\left({
      \frac{\sqrt{\del}}{2\pi}+\frac{|y|}{2\sqrt{2\pi}}
      +\frac{\sqrt{\del}}{4\pi}(x+y)^2
      +\frac{y^2}{4\pi\sqrt{\del}}
    }\right)\frac{70}{\sqrt{2\pi}}|x|\\
  &\quad+10\left[{
      (|x|+|y|)^3+\frac{|x||y|^2}{\sqrt{\del}}(|x|+|y|)
      +\frac{|y|^4}{\del^{3/2}}+\del^{3/2}(x+y)^2+\del(|x|+|y|)
    }\right]
  \end{split}\end{equation*}
by \eqref{65d}.

Now suppose that $\del\le 1$ and $\al,\beta\in\RR$. Let
$y=\del^{1/2+\al}$, $x=-\del^{1/4+\beta}$, and assume that %
$y\le-x\le 1$. Using the fact that $|x|+|y|\le 2|x|\le 2$, we have
  \begin{align*}
  |R_\del(x,y)|&\le\frac{50}{\sqrt{2\pi}}|x|^3
    +\frac{70}{(2\pi)^{3/2}}|x|\left({
      \sqrt{\del}+2|y|+2\sqrt{\del}+\frac{y^2}{\sqrt{\del}}
    }\right)\\
  &\quad+10\left({
      8|x|^3+2\frac{|x||y|^2}{\sqrt{\del}}
      +\frac{|y|^4}{\del^{3/2}}+4\del^{3/2}x^2+2\del|x|
    }\right)\\
  &\le 25|x|^3+5|x|\left({
      3\sqrt{\del}+2|y|+\frac{y^2}{\sqrt{\del}}
    }\right)\\
  &\quad+80|x|^3+20\frac{|x||y|^2}{\sqrt{\del}}
    +10\frac{|y|^4}{\del^{3/2}}+40\del^{3/2}x^2+20\del|x|
  \end{align*}
which reduces to
  \begin{equation*}\begin{split}
  |R_\del(x,y)|&\le 105\del^{3/4+3\beta}+15\del^{3/4+\beta}
    +10\del^{3/4+\al+\beta}+25\del^{3/4+2\al+\beta}\\
  &\quad+10\del^{1/2+4\al}+40\del^{2+2\beta}
    +20\del^{5/4+\beta}.
  \end{split}\end{equation*}
To simplify further, suppose $\al>0$. Then
  \begin{align*}
  |R_\del(x,y)|&\le 105\del^{3/4+3\beta}+15\del^{3/4+\beta}
    +10\delta^{3/4+\beta}+25\del^{3/4+\beta}\\
  &\quad+10\del^{1/2+4\al}+40\del^{2+2\beta}
    +20\del^{3/4+\beta}\\
  &= 105\del^{3/4+3\beta}+70\del^{3/4+\beta}
    +10\del^{1/2+4\al}+40\del^{2+2\beta}.
  \end{align*}
Now, if $\beta\ge 0$, then $2+2\beta>3/4+\beta$, and
$|R_\del(x,y)|\le
115\del^{3/4+3\beta}+110\del^{3/4+\beta}+10\del^{1/2+4\al}$.
Otherwise, if $\beta<0$, then $2+2\beta>3/4+3\beta$, and
$|R_\del(x,y)|\le
145\del^{3/4+3\beta}+70\del^{3/4+\beta}+10\del^{1/2+4\al}$. In
either case,
  \[
  |R_\del(x,y)|\le
    150(\del^{3/4+3\beta}+\del^{3/4+\beta}+\del^{1/2+4\al})
  \]
whenever $\al>0$. On the other hand, suppose $\al<0$. Then
  \begin{align*}
  |R_\del(x,y)|&\le 105\del^{3/4+3\beta}+15\del^{3/4+2\al+\beta}
    +10\delta^{3/4+2\al+\beta}+25\del^{3/4+2\al+\beta}\\
  &\quad+10\del^{1/2+4\al}+40\del^{2+2\beta}
    +20\del^{3/4+2\al+\beta}\\
  &= 105\del^{3/4+3\beta}+70\del^{3/4+2\al+\beta}
    +10\del^{1/2+4\al}+40\del^{2+2\beta}.
  \end{align*}
If $\beta\ge 0$, then $2+2\beta>3/4+\beta\ge 3/4+2\al+\beta$; if
$\beta<0$, then $2+2\beta>3/4+3\beta$. We therefore have
  \[
  |R_\del(x,y)|\le 150(\del^{3/4+3\beta}+\del^{3/4+2\al+\beta}
    +\del^{1/2+4\al})
  \]
whenever $\al<0$.

In summary, we have an expansion for $p_1(x,y,\del)$ given by
\eqref{70a}, together with a remainder estimate of the form
  \begin{equation}
  |R_\del(x,y)|\le 150(\del^{3/4+3\beta}
    +\del^{3/4+2(\al\wedge 0)+\beta}+\del^{1/2+4\al}),\label{73a}
  \end{equation}
valid for $0<\del\le 1$ whenever $y=\del^{1/2+\al}$ and
$x=-\del^{1/4+\beta}$ satisfy $y\le-x\le 1$. Moreover, by
symmetry, the same bound holds for $|R_\del(-x,-y)|$.

Now fix $\Del>0$. Choose $\del_0\le 1$ such that
  \begin{equation}
  900(\del_0^{1/4}\vee\del_0^{3\Del})<(2\pi)^{-1/2}.\label{74a}
  \end{equation}
Let $\al\ge\Del$ and $0<\del\le\del_0$. Set $\beta=\al$,
$y=\del^{1/2+\al}$, and $x=-\del^{1/4+\beta}$. Note that by
\eqref{37a}-\eqref{38d}
  \[
  \mut(x,y,\del)=p_1(x,y,\del)-p_1(-x,-y,\del),
  \]
so by \eqref{70a}
  \[
  \mut=\frac{2y}{\sqrt{2\pi}}+R_\del(x,y)-R_\del(-x,-y).
  \]
Since $\del\le 1$, we have $y\le-x\le 1$. Hence, by \eqref{73a}
and \eqref{74a},
  \begin{align*}
  |R_\del(x,y)-R_\del(-x,-y)|
    &\le 300(2\del^{3/4+\al}+\del^{1/2+4\al})\\
  &= 300(2\del^{1/4}+\del^{3\al})y\\
  &\le 300(2\del_0^{1/4}+\del_0^{3\Del})y\\
  &\le 900(\del_0^{1/4}\vee\del_0^{3\Del})y<(2\pi)^{-1/2}y.
  \end{align*}
Therefore, $\mut\ge(2\pi)^{-1/2}y$, which proves (i).

For (ii), observe that $\mut>0$ implies $q_1<q_2$. Hence
$\eet=p_1q_2+p_2q_1\le 2q_2$. Moreover,
  \begin{equation}\label{76a}
  \begin{split}
  q_2&=1-p_1(-x,-y,\del)\\
  &\le\frac1{\pi}\tan^{-1}\sqrt{\del}+\frac{|y|}{\sqrt{2\pi}}
    +\frac{\sqrt{\del}}{2\pi}(x+y)^2+\frac{y^2}{2\pi\sqrt{\del}}
    +|R_\del(-x,-y)|\\
  &\le\del^{1/2}+\del^{1/2+\al}+\del^{1+2\al}+\del^{1/2+2\al}
    +150(\del^{3/4+3\al}+\del^{3/4+\al}+\del^{1/2+4\al})\\
  &\le 500\,\del^{1/2},
  \end{split}
  \end{equation}
so $\eet\le 1000\,\del^{1/2}$. By making $\del_0$ smaller if
necessary we can ensure that $1000\,\del^{1/2}<1/2$. \qed

\begin{lemma}\label{L:small}
Let $p>2$. Fix $0<\Del<1/2$ and let $\del_0$ be as in Lemma
\ref{L:taylor}. Suppose $\ee>0$, $0<\del\le\del_0$, and $n\ge 3$
satisfy $\ee/\sqrt{n}\le\del^{1/2+\Del}$. Then
  \[
  P\left({M_n(1+\del)-M_n(1)>\frac{\ee}{\sqrt{n}}}\right)
    \le C(\ee^{-1}\del^{1/4})^p,
  \]
where $C$ depends only on $p$ and $\Del$.
\end{lemma}

\pf Let $y=\ee/\sqrt{n}$ and choose $\al\ge\Del$ such that
$y=\del^{1/2+\al}$. Set $x_0=-\del^{1/4+\al}$. By Corollary
\ref{C:cond}, Lemma \ref{L:balance}, Lemma \ref{L:chebyplus}, and
Lemma \ref{L:taylor},
  \[
  P\left({M_n(1+\del)-M_n(1)>\frac{\ee}{\sqrt{n}}}\right)
    \le C_{p/2}\frac{\eet^{p/2}}{(k-1)^{p/2}\mut^p}
    +C_p(\ee^{-1}\del^{1/4})^p,
  \]
where $\eet=\eet(x_0,y,\del)\le 1000\,\del^{1/2}<1/2$ and
  \[
  \mut=\mut(x_0,y,\del)\ge\frac1{\sqrt{2\pi}}\,\del^{1/2+\al}
    =\frac1{\sqrt{2\pi}}\cdot\frac{\ee}{\sqrt{n}}>0.
  \]
Hence,
  \[
  \frac{\eet^{p/2}}{(k-1)^{p/2}\mut^p}
    \le C\frac{\del^{p/4}}{n^{p/2}(\ee/\sqrt{n})^p}
    =C(\ee^{-1}\del^{1/4})^p,
  \]
which completes the proof. \qed

\section{The Medium Jump Regime and Final Proof}

Our analysis of the medium jump regime will require only minor
modifications to the methods of Section \ref{small}.

\begin{lemma}\label{L:taylormod}
Fix $0<\Del<1/16$ and set $\Del'=(1-16\Del)/12>0$. Then there
exists $\del_0>0$ such that
  \begin{description}
  \item{(i)} $\mut(-\del^{1/4+\al},\del^{1/2+\al},\del)
    \ge\frac1{\sqrt{2\pi}}\,\del^{1/2+\Del}$, and
  \item{(ii)} $\eet(-\del^{1/4+\al},\del^{1/2+\al},\del)
    \le 1000\,\del^{1/2-4\Del'}<\frac12$
  \end{description}
for all $-\Del'\le\al\le\Del$ and all $0<\del\le\del_0$.
\end{lemma}

\pf For fixed $0<\Del<1/16$, choose $\del_0>0$ as in Lemma
\ref{L:taylor}. By \eqref{37a}, $p_1$ is increasing in $y$. Hence,
if $x=-\del^{1/4+\al}$ and $y=\del^{1/2+\Del}$, then by
\eqref{70a} and \eqref{73a},
  \begin{align*}
  \mut(x,\del^{1/2+\al},\del)
    &=p_1(x,\del^{1/2+\al},\del)-p_1(-x,-\del^{1/2+\al},\del)\\
  &\ge p_1(x,y,\del)-p_1(-x,-y,\del)\\
  &=\frac{2y}{\sqrt{2\pi}}+R_\del(x,y)-R_\del(-x,-y),
  \end{align*}
where
  \begin{align*}
  |R_\del(x,y)-R_\del(-x,-y)|
    &\le 300(\del^{3/4+3\al}+\del^{3/4+\al}+\del^{1/2+4\Del})\\
  &\le 300(2\del^{3/4-3\Del'}+\del^{1/2+4\Del}).
  \end{align*}
However, note that $3/4-3\Del'=1/2+4\Del$. Hence, by \eqref{74a},
  \[
  |R_\del(x,y)-R_\del(-x,-y)|\le 900\del^{1/2+4\Del}
    =900\del^{3\Del}y<(2\pi)^{-1/2}y.
  \]
Therefore, $\mut\ge(2\pi)^{-1/2}y$, which proves (i).

For (ii), observe that $\eet\le 2q_2$ and, as in \eqref{76a},
  \begin{align*}
  q_2&=1-p_1(-x,-\del^{1/2+\al},\del)\\
  &\le\del^{1/2}+\del^{1/2+\al}+\del^{1+2\al}+\del^{1/2+2\al}
    +150(\del^{3/4+3\al}+\del^{3/4+\al}+\del^{1/2+4\al})\\
  &\le 4\del^{1/2-2\Del'}
    +150(2\del^{3/4-3\Del'}+\del^{1/2-4\Del'})\\
  &\le 500\,\del^{1/2-4\Del'}.
  \end{align*}
Note that $1/2-4\Del'>1/6$, so that by making $\del_0$ smaller if
necessary, we can ensure that $1000\,\del^{1/2-4\Del'}<1/2$. \qed

\begin{lemma}\label{L:medium}
Fix $p>2$. Let $\Del=1/18$ and choose $\del_0>0$ as in Lemma
\ref{L:taylormod}. Suppose $\ee>0$, $0<\del\le\del_0$, and $n\ge
3$ satisfy $\del^{5/9}\le\ee/\sqrt{n}\le\del^{53/108}$. Then
  \[
  P\left({M_n(1+\del)-M_n(1)>\frac{\ee}{\sqrt{n}}}\right)
    \le C(\ee^{-1}\del^{1/6})^p,
  \]
where $C$ depends only on $p$.
\end{lemma}

\pf Let $\Del=1/18$ and $\Del'=(1-16\Del)/12=1/108$ and observe
that $y=\ee/\sqrt{n}=\del^{1/2+\al}$ for some
$\al\in[-\Del',\Del]$. Set $x_0=-\del^{1/4+\al}$. By Corollary
\ref{C:cond}, Lemma \ref{L:balance}, Lemma \ref{L:chebyplus}, and
Lemma \ref{L:taylormod},
\[
  P\left({M_n(1+\del)-M_n(1)>\frac{\ee}{\sqrt{n}}}\right)
    \le C_{p/2}\frac{\eet^{p/2}}{(k-1)^{p/2}\mut^p}
    +C_p(\ee^{-1}\del^{1/4})^p,
  \]
where $\eet=\eet(x_0,y,\del)\le 1000\,\del^{1/2-4\Del'}<1/2$ and
  \[
  \mut=\mut(x_0,y,\del)
    \ge\frac1{\sqrt{2\pi}}\,\del^{1/2+\Del}>0.
  \]
Note that $n=\ee^2y^{-2}=\ee^2\del^{-1-2\al}$. Hence,
  \begin{align*}
  \frac{\eet^{p/2}}{(k-1)^{p/2}\mut^p}
    &\le C(\eet\ee^{-2}\del^{1+2\al}\mut^{-2})^{p/2}\\
  &\le C(\del^{1/2-4\Del'}\ee^{-2}\del^{1-2\Del'}
    \del^{-1-2\Del})^{p/2}\\
  &=C(\ee^{-2}\del^{1/2-6\Del'-2\Del})^{p/2}.
  \end{align*}
Since $1/2-6\Del'-2\Del=1/3$, this completes the proof. \qed

With the completion of our lemmas, we have made short work of the
only proof that remains.

\bigskip\noindent{\bf Proof of Lemma \ref{L:key}:} Take
$\Del=1/108$ in Lemma \ref{L:large} and, for each $p>2$, let
$C_{p,1}$ be the constant that appears in that lemma. Then take
$\Del=1/18$ in Lemma \ref{L:taylormod}. Let $\del_0>0$ be as in
that lemma and note that the conclusions of Lemmas \ref{L:small}
and \ref{L:medium} hold for this choice of $\del_0$. For each
$p>2$, let $C_{p,2}$ be the larger of the constants appearing in
those two lemmas and let $C_p=C_{p,1}\vee C_{p,2}$.

Now let $0<\ee<1$, $0<\del\le\del_0$, and $n\ge 3$. Choose
$\al>-1/2$ such that $\ee/\sqrt{n}=\del^{1/2+\al}$. If
$\al\le-1/108$, then by Lemma \ref{L:large},
  \[
  P\left({M_n(1+\del)-M_n(1)>\frac{\ee}{\sqrt{n}}}\right)
    \le C_{p,1}(\ee^{-1}\del^{1/4})^p
    \le C_p(\ee^{-1}\del^{1/6})^p.
  \]
If $\al\ge 1/18$, then by Lemma \ref{L:small},
  \[
  P\left({M_n(1+\del)-M_n(1)>\frac{\ee}{\sqrt{n}}}\right)
    \le C_{p,2}(\ee^{-1}\del^{1/4})^p
    \le C_p(\ee^{-1}\del^{1/6})^p.
  \]
If $-1/108\le\al\le 1/18$, then by Lemma \ref{L:medium},
  \[
  P\left({M_n(1+\del)-M_n(1)>\frac{\ee}{\sqrt{n}}}\right)
    \le C_{p,2}(\ee^{-1}\del^{1/6})^p
    \le C_p(\ee^{-1}\del^{1/6})^p,
  \]
and we are done. \qed

\bigskip\noindent{\bf Acknowledgments.} This material appeared in
my doctoral dissertation and I would like to thank my advisors,
Chris Burdzy and Zhen-Qing Chen, for their instruction and
guidance. I also thank the referees for their many helpful
suggestions. For helpful discussions and for sharing their
insights, I would like to thank Jon Wellner, Davar Khoshnevisan,
Wenbo Li, Bruce Erickson, Yaozhong Hu, and especially Tom Kurtz,
to whom I am indebted for his assistance in preparing this article
for publication. This work was done while supported by VIGRE, for
which I thank the NSF, University of Washington, and University of
Wisconsin-Madison.


\begin{thebibliography}{9}

\bibitem{D} Detlef D\"urr, Sheldon Goldstein, and Joel L.
Lebowitz, Asymptotics of Particle Trajectories in Infinite
One-Dimensional Systems with Collisions. {\it Communications on
Pure and Applied Mathematics}, {\bf 38} (1985), 575--597.

\bibitem{H} T. E. Harris, Diffusions with collisions between
particles. {\it Journal of Applied Probability}, {\bf 2} (1965),
323--338.

\bibitem{K} Ioannis Karatzas and Steven E. Shreve, {\it Brownian
Motion and Stochastic Calculus}. Springer, 1991

\bibitem{R} R. D. Reiss, {\it Approximate Distributions of Order
Statistics With Applications to Nonparametric Statistics}.
Springer, 1989

\bibitem{S} F. Spitzer, Uniform motion with elastic collisions of
an infinite particle system. {\it J. Math. Mech.}, {\bf 18}
(1968), 973--989.

\bibitem{St} Daniel W. Stroock, {\it Probability Theory, An
Analytic View}. Cambridge University Press, 1993

\bibitem{Sw} Jason Swanson, {\it Variations of Stochastic
Processes: Alternative Approaches}. Doctoral Dissertation,
University of Washington, 2004.

\bibitem{T} P. F. Tupper, A test problem for molecular dynamics
integrators. {\it IMA Journal of Numerical Analysis}, {\bf 25(2)}
(2005), 286--309.

\end{thebibliography}
\end{document}